\documentclass[12pt,leqno]{amsart}
\usepackage{amsfonts,amsmath,amsthm,amssymb,bbm,eufrak}
\usepackage[usenames,dvipsnames]{color}
\usepackage{enumerate,epsfig,graphics,indentfirst}
\usepackage{hyperref,latexsym,multicol}
\usepackage{pxfonts,rotating,yfonts,ulem}
\usepackage[labelsep=space]{caption}
\usepackage[all]{xy}
\usepackage{tikz}
\usetikzlibrary{decorations.markings}
\usepackage[labelformat=empty]{caption}
\usepackage{mathtools}

\newtheorem{theorem}{Theorem}[section]
\newtheorem{lemma}[theorem]{Lemma}
\newtheorem{proposition}[theorem]{Proposition}
\newtheorem{corollary}[theorem]{Corollary}

\theoremstyle{definition}

\usepackage[ansinew]{inputenc}
\usepackage{ifthen}
\usepackage{animate}

\newcommand{\be}{\begin{equation}}
\newcommand{\ee}{\end{equation}}

\begin{document}

\newpage

\title[Estimates for Thue inequalities with few coefficients]{Explicit and Mixed Estimates for Thue inequalities with few coefficients}
\author[Saradha]{N. Saradha}
\address[1]{Fellow, Indian Academy of Sciences\\
B-706, Everard Towers\\
Eastern Express Highway\\
Sion, Mumbai 400\,022 \textsc{India}}
\email{saradha54@gmail.com}

\author[Sharma]{Divyum Sharma}
\address[2]{Department of Mathematics\\
Birla~Inst{i}tute~of~Technology~and~Science, Pilani 333\,031 \textsc{India}}
\email{\texttt{divyum.sharma\symbol{64}pilani.bits-pilani.ac.in}}

\date{}

\begin{abstract}
Let $F(x,y)$ be an irreducible form of degree $r\geq 3$ and having $s+1$ non-zero coefficients. Let $h\geq 1$ be an integer and consider the Thue inequality
$$|F(x,y)|\leq h.$$
Following the seminal work of Thue in 1909, several papers were written giving an upper bound for the number of solutions of the above inequality as $\ll c(r,s,h)$ where $c(r,s,h)$ is an explicit function of $r,s$ and $h.$ Invariably, the absolute constant involved in $\ll$ has been left undetermined. In  this paper, following Bombieri, Schmidt and Mueller, we give three different upper bounds which are explicit in every aspect.
\end{abstract}

\maketitle
\section{Introduction} 
\vskip 2mm
\noindent
{\bf A brief history :}
Let $F(x,y)=a_0x^r+a_1x^{r-1}y+\cdots+a_ry^r\in \mathbb Z[x,y]$ be an irreducible homogeneous polynomial of degree $r\geq 3.$  Let $s+1$ be the number of non-zero $a_i$'s, $D(F)$ the discriminant of $F$, $H(F)$ the na\"ive height of $F$  and $M(F)$ the Mahler measure of $F.$ 
For any positive integer $h \geq 1,$ we consider the inequality
\be\label{inequality}
|F(x,y)|\leq h.
\ee
We want to estimate the number of solutions $N(F,h)$ of \eqref{inequality}. 
Siegel \cite{Sie} proposed that $N(F,h)$ can be bounded by a function depending only on $r$ and $h,$ otherwise independent of $F.$ 
This has been the theme of several papers since the seminal work of Thue \cite{Thue} in 1909. See in particular, the papers of Evertse \cite{Ever}, Bombieri \& Schmidt 
\cite{Bomb-Sch},
Schmidt \cite{Sch} and Mueller \& Schmidt \cite{Mue-Sch}. The first two papers, in fact deal with estimating the number of solutions when the inequality is replaced by equality in \eqref{inequality}. It is not difficult to modify their proofs to deal with \eqref{inequality}.
 In \cite{Sch}
and \cite{Mue-Sch},
it was shown that $N(F,h)$ satisfies
\be\label{rs}
N(F,h)\ll (rs)^{1/2}h^{2/r}(1+\log h^{1/r})
\ee
and
\be\label{ss}
N(F,h)\ll s^2h^{2/r}(1+\log h^{1/r}),
\ee
respectively. Here, the constants implied by $\ll$ are absolute but not explicitly given.
While dealing with inequality \eqref{inequality}, the term $h^{2/r}$ in the above bounds is needed. The logarithmic factor was improved by Thunder \cite{Thun} as
$$
N(F,h)\ll (rs)^{1/2}\left(h^{1/r}+h^{2/r}|D(F)|^{-1/r(r-1)}\right)\bigg(1+\frac{\log\log h}{\log r}\bigg).
$$
In recent times there are some papers improving $s^2$ in \eqref{ss}. See
\cite{Akh-pam},\cite{Ben},\cite{Sara-Div1} and \cite{Sara-Div2}. 
\vskip 2mm
\noindent
{\bf New Explicit Estimates :}
In all the above results the absolute constants implied by $\ll$ are not explicitly given. There are several instances when such constants are explicitly computed, especially when $|D(F)|$ is large. See for instance, \cite{Bomb-Sch},\cite{Gyor} and \cite{Div-Sar}.
In 1970, Siegel \cite{Sie} showed that if $F$ is a {\it diagonalizable} form and 
$$
|D(F)|\geq 2^{r^2-r}r^{r+4(45+593/913)}h^{2(r-1)+(45+593/913)r}
$$
 then
$$N(F,h)\leq 2rh^{2/r}(1+\log h^{1/r}).$$
The condition on $|D(F)|$ was improved in \cite{Akh-Sara-Div} where  precise statements
can be found. See also \cite{Sara-Div3}. Then one wonders how big the absolute constants 
would be in \eqref{rs} and \eqref{ss}. For any real number $x,$ put
$$\log^*x=
\begin{cases}
1\ \textrm{if}\ x\leq e,\\
\log x\ \textrm{if} \ x> e.
\end{cases}
$$
Let 
\be\label{zr}
z(r)=68\log^3r+150\log^2r+\log r.
\ee
We show
{\color{black}\begin{theorem}\label{Main1}
Suppose $F$ has at most $s+1$ non-zero terms. 
Then
\begin{align*}
 N(F,h)&\leq 2(1+\log^*h^{1/r})\bigg(3+4h^{1/r}+\frac{e^{105}h^{2/r}}{|D(F)|^{1/r(r-1)}}\bigg)\times\\
 &\left(1+8r+48s\bigg(15+\frac{1}{\log(r-1)}\log\log^*(2h^{1/r})\bigg)\right).
\end{align*}
\end{theorem}}
In the next two theorems we relax the dependency on $r.$ 

{\color{black}\begin{theorem}\label{Main2}
Suppose $F$ has at most $s+1$ non-zero terms. Assume that
$r\geq 132434.$
Then
\begin{align*}
&N(F,h)\leq \\
&2(1+\log^*h^{1/r})\bigg(3+4h^{1/r}+\frac{e^{105}h^{2/r}}{|D(F)|^{1/r(r-1)}}\bigg)\times\\
&\bigg(1+40\big(6+z(r)+25(sr)^{1/2}\big)+48s\big(15+\frac{1}{\log(r-1)}\log\log^*(2h^{1/r})\big)\bigg).
\end{align*}
\end{theorem}}

Define
$$
\epsilon=
\begin{cases}

0\ \textrm {if} \ r \geq 3s^3,\\
1\ otherwise.
\end{cases}
$$
{\color{black}\begin{theorem}\label{Main3}
Suppose $F$ has at most $s+1$ non-zero terms with $r\geq 3s.$ Then
\begin{align*}
N(F,h)&\leq   2\left(1+\log^* h^{1/r}\right) h^{1/r}\times
\bigg(1+324s+648s (\log^*s)^\epsilon+\\
&16z(r)+16\log^* (2h^{1/r})+36s\log\bigg(1+\frac{\log^* (2^{1/(r-s)}h^{1/r(r-s)})}{\log^* H(F)}\bigg)\bigg).
\end{align*}
\end{theorem}}
\section{Counting solutions of \eqref{inequality}}
Observe that $N(F,h)$ is an increasing function of $h.$ The content of $F$, denoted as $c(F)$, is $\gcd(a_0,\ldots,a_r).$
Suppose $c(F)>1.$ Then write
$$F(x,y)=c(F)G(x,y)$$
where $c(G)=1.$ Then we may consider \eqref{inequality} with $h$ replaced by $c(F)^{-1}h$ and $F$ by $G.$ Thus $N(F,h)=N(G,c(F)^{-1}h)\leq N(G,h).$ Hence it is no loss of generality to assume that $c(F)=1.$

\vskip 2mm
\noindent
{\bf A reduction of the problem :}

Let us denote by $P(F,h)$ the number of primitive solutions i.e., solutions $(x,y)$ of \eqref{inequality} satisfying $\gcd(x,y)=1.$ We put $P(F,0)=0.$ Obviously $N(F,h)\geq P(F,h).$ 
Following the argument of Schmidt \cite[p. 245-246]{Sch}
we have the ensuing lemma.
\begin{lemma}\label{PFh}
We have
$$ (i) \ N(F,h)\leq P(F,h) h^{1/r}.$$
More precisely, we have
$$(ii)\ N(F,h)\leq  P(F,h)+h^{1/r}\sum_{m=1}^{h-1}P(F,m)r^{-1}m^{-1-1/r}.$$
\end{lemma}
\begin{proof}
Let $\pi(F,h)$ denote the number of primitive solutions of 
$$|F(x,y)|=h.$$
Then 
\begin{align}
N(F,h)&=\sum_{m=1}^h\pi(F,m)\bigg\lfloor{\bigg(\frac{h}{m}\bigg)^{1/r}\bigg\rfloor}\nonumber\\
&\leq h^{1/r}\sum_{m=1}^h\pi(F,m)m^{-1/r}\nonumber\\
&=h^{1/r}\sum_{m=1}^h\big(P(F,m)-P(F,m-1)\big){m}^{-1/r}\nonumber\\
&=P(F,h)+h^{1/r}\sum_{m=1}^{h-1}P(F,m)\big(m^{-1/r}-(m+1)^{-1/r}\big)\label{mm1}\\
&\leq P(F,h)+h^{1/r}P(F,h)\sum_{m=1}^{h-1}\big(m^{-1/r}-(m+1)^{-1/r}\big)\nonumber\\
&\leq P(F,h)+h^{1/r}P(F,h)\bigg(1-h^{-1/r}\bigg)=P(F,h)h^{1/r}.\nonumber
\end{align}

We can be more precise. By the mean value theorem, there exists $\theta$ with $m<\theta<m+1$ such that
$$m^{-1/r}-(m+1)^{-1/r}=r^{-1}\theta^{-1-1/r}< r^{-1} m^{-1-1/r}.$$
Thus by \eqref{mm1},
\be\label{NPP}
N(F,h)\leq P(F,h)+h^{1/r}\sum_{m=1}^{h-1}P(F,m)r^{-1}m^{-1-1/r},
\ee
which proves (ii). 
\end{proof}
The above lemma enables us to focus on estimating the primitive solutions $P(F,h)$ of \eqref{inequality}.  
Denote by $P'(F,h)$ the number of primitive solutions of the inequality
\be\label{secondinequality}
\frac{h}{2^r}\leq |F(x,y)|<h.
\ee
For any positive integer $m,$ let $u(m)$ be the unique integer defined by
\be\label{um}
2^{ru(m)}\leq m<2^{r(u(m)+1)}.
\ee
Then  $u(m)$ is a non-decreasing function of $m$.
Observe that if $h=2^{rj},$ then $u(h)=j.$ Note that
\be\label{PP}
P(F,h)\leq \sum_{j=0}^{u(h)+1}P'(F,2^{rj}).
\ee
So the problem is now reduced to estimating $P'(F,m)$ for any integer $m\geq 1.$ 

\section{Action of $M_2(\mathbb Z)^{\times}$ on Forms}
Let $\alpha_1,\ldots,\alpha_r$ be the roots of the polynomial
$$f(x)=F(x,1)=a_0x^r+a_1x^{r-1}+\cdots+a_r$$
then we have 
$$F(x,y)=a_0(x-\alpha_1y)\cdots(x-\alpha_ry).$$
The discriminant $D(F)$ of the polynomial $F(x,1)$ is given by
$$D(F)=a_0^{2r-2}\prod_{1\leq i<j\leq r}(\alpha_i-\alpha_j)^2.$$
Let $M_2(\mathbb Z)^{\times}$ denote the set of all $2\times 2$ integer matrices with non-zero determinant.
For 
$A=\begin{pmatrix}{a\ \ b}\\{c\ \ d}\end{pmatrix}\in M_2(\mathbb Z)^{\times},$ set $F_A(x,y)$ as
$$F_A(x,y)=F(ax+by,cx+dy).$$
Then 
\be\label{DFA}
D(F_A)=\det(A)^{r(r-1)}D(F).
\ee
For $p\geq 1$, let
\be\label{Aj}
A_0=\begin{pmatrix}{p\ \ 0}\\{0\ \ 1}\end{pmatrix},
A_j=\begin{pmatrix}{0\ \ -1}\\{p\ \ \ \ \ j}\end{pmatrix}, j=1,\ldots, p.
\ee
The following decomposition of $\mathbb Z^2$ is well known.
\begin{lemma}\label{decomp} Let $p$ be any prime. Then
$$\mathbb Z^2=\cup_{j=0}^pA_j\mathbb Z^2$$
where $A_j, 0\leq j\leq p$ are matrices given in \eqref{Aj}.
\end{lemma}
Now suppose $(x,y)$ is a solution of \eqref{secondinequality} and 
$(x,y)=A_j(x',y').$ Then
$$|F_{A_j}(x',y')|=|F(A_j(x',y'))|=|F(x,y)|.$$
Also, $|D(F_{A_j})|=p^{r(r-1)}|D(F)|.$ 
Thus if $n_j$ denotes the number of primitive solutions of  \eqref{secondinequality} with $F$ replaced by $F_{A_j},$ then
\be\label{pinequality}
P'(F,m)\leq n_0+\cdots+n_p\leq (p+1)\max_{0\leq j\leq p} P'(F_{A_j},m).
\ee
This was shown in \cite{Bomb-Sch}.
A nice property of $F_{A_j}$ is that it has discriminant of absolute value at least $p^{r(r-1)}.$
Thus the problem reduces to judicious choice of $p$ and then estimating $P'(G,m)$ where
$G$ is a form having large discriminant. 

Let $C$ be a class of forms which is closed under the action of $M_2(\mathbb Z)^{\times}.$
For any given integer $p,$ let $C_p$ be the subclass of forms $F\in C$ with  
$$|D(F)|\geq p^{r(r-1)}.$$
By \eqref{DFA}, $C_p$ is closed under the action of $M_2(\mathbb Z)^{\times}.$ Note that for $F\in C,$ we have $F_{A_j}\in C_p.$

Let us consider forms with at most $s+1$ non-zero coefficients. The $M_2(\mathbb Z)^{\times}$ actions do not preserve the number of non-zero coefficients of $F.$ In order to circumvent this, Schmidt introduced a class $C(t)$ as follows.
Let $C(t)$ denote a class of irreducible forms $F(x,y)\in \mathbb Z[x,y]$ of degree $r$ and such that the form
$$uF_x+vF_y$$
has at most $t$ real zeros for real numbers $u$ and $v$ with $(u,v)\neq (0,0).$
Here $F_x$ and $F_y$ denote partial derivatives with respect to $x$ and $y$ respectively.
Note that $uF_x+vF_y$ is a form of degree $r-1$ and is not identically zero, since $F$ is irreducible.   Suppose $G=F_A,$ then
\begin{align}
&uG_x(x,y)+vG_y(x,y)\label{uv}\\
&=(ua+vb)F_x(ax+by,cx+dy)+(uc+vd)F_y(ax+by,cx+dy)\nonumber\\
&=u_1F_x(ax+by,cx+dy)+v_1F_y(ax+by,cx+dy)\nonumber
\end{align}
with $u_1=ua+vb,v_1=uc+vd.$ Thus $C(t)$ is closed under $M_2(\mathbb Z)^{\times}.$
The following is \cite[Lemma 2]{Sch}.
\begin{lemma}\label{4s}
Suppose $F(x,y)\in \mathbb Z[x,y]$ is irreducible of degree $r$ and having at most $s+1$ non-vanishing coefficients. Then $F(x,y)\in C(4s-2).$
\end{lemma}

Let $f(x)=F(x,1).$ Since $F(x,y)$ is irreducible, we have $f(0)\neq 0.$ Hence by \cite[Lemma 1]{Sch}, $f(x)$ has no more than $2s$ real roots, and $f'(x) $ has no more than $2s-1$ real roots.
So $ff'$ has $\leq 4s-1$ real roots.
Suppose $F(x,y)\in C(t).$ Then $f'(x)=1. F_x(x,1)$ has no more than $t$ real roots. This implies then that $f(x)$ has at most $t+1$ real roots.
So together with \eqref{uv}, the following corollary is immediate.
\begin{corollary}\label{realroots}
Suppose $F(x,y)\in \mathbb Z[x,y]$ is irreducible of degree $r.$ Let $f(x)=F(x,1).$ Then 

 (i) $ff'$ has at most $4s-1$ real roots if $F(x,y)$ has at most $s+1$ non-zero coefficients.

(ii) $ff'$ has at most $2t+1$ real roots if $F(x,y)\in C(t).$

(iii) $F_{A_j}\in C_p(t)$ for all $j=0,1,\ldots,p$ if $F(x,y)\in C(t).$
\end{corollary}
{\bf Choice of the prime $p$ :}
Let $q(r)$ be any function of $r$. Given $h\geq 1,$ choose the least prime $p$ exceeding 
\be\label{uo}
U_0=\frac{e^{q(r)}2^{2u(h)-1}}{|D(F)|^{1/r(r-1)}}.
\ee
We will take a specific function $q(r)$ later.
If $U_0< 2,$ then $p=2.$ If $U_0\geq 2,$ then
\be\label{lbp}
p> \frac{e^{q(r)}h^{2/r}}{8|D(F)|^{1/r(r-1)}}.
\ee
Put
\be\label{u}
U_1=3+\frac{e^{q(r)}h^{2/r}}{|D(F)|^{1/r(r-1)}}.
\ee
If $U_0\geq 2,$ then by Bertrand's postulate, there exists a prime $p'$ such that
$$U_0<p'<2U_0.$$
Hence
$$
p+1\leq p'+1<2U_0+1=\frac{e^{q(r)}2^{2u(h)}}{|D(F)|^{1/r(r-1)}}+1\leq 
\frac{e^{q(r)}h^{2/r}}{|D(F)|^{1/r(r-1)}}+1.
$$
Thus we have $p+1\leq U_1$ always.
This gives by \eqref{pinequality} that
\be\label{Ugm}
P'(F,h)\leq U_1\max_{0\leq j\leq p} P'(F_{A_j},h).
\ee

Combining this with \eqref{PP} and Lemma \eqref{PFh}(ii) we have the following lemma.
\begin{lemma}\label{NCUg}
Let $h\geq 1$ and 
let $U_0$ and $U_1$ be given by \eqref{uo} and \eqref{u}, respectively. Suppose $p$ is the least prime exceeding $U_0.$ Assume that
\be\label{Cgm}
\max_{0\leq j\leq p}P'(F_{A_j},h)\leq g(h)
\ee
where $g(h)$ is an increasing function of $h.$ Then
$$
P(F,h)\leq 2g(2^rh)(1+\log^* h^{1/r})\left(3+\frac{e^{q(r)}h^{2/r}}{|D(F)|^{1/r(r-1)}}\right).
$$
and
$$
N(F,h) \leq  2g(2^rh) (1+\log^* h^{1/r}) \left(3+4h^{1/r}+ \frac{2e^{q(r)}h^{2/r}}{|D(F)|^{1/r(r-1)}}\right).
$$
\end{lemma}
\begin{proof}
Using \eqref{Ugm} and \eqref{Cgm}, we calculate 
\begin{align*}
\sum_{j=0}^{u(h)+1} P'(F,2^{rj})\leq & U_1\sum_{j=0}^{u(h)+1} g(2^{rj})\\
\leq &\left(3+\frac{e^{q(r)}h^{2/r}}{|D(F)|^{1/r(r-1)}}\right) g(2^{r(u(h)+1)})\bigg(u(h)+2\bigg)\\
\leq & \left(3+\frac{e^{q(r)}h^{2/r}}{|D(F)|^{1/r(r-1)}}\right)g(2^rh)\bigg( 2+\frac{\log h}{r\log 2}\bigg).
\end{align*}
Thus by \eqref{PP},
\be\label{63}
P(F,h)\leq 2\left(3+\frac{e^{q(r)}h^{2/r}}{|D(F)|^{1/r(r-1)}}\right)g(2^rh)\bigg( 1+\log^* h^{1/r}\bigg)
\ee
which proves the bound for $P(F,h).$

To find a bound for $N(F,h)$ we use Lemma \ref{PFh} (ii).  First note that
$$\sum_{m=1}^{h-1} \frac{1}{m^{1+1/r}}\leq 1+\int_1^\infty \frac{1}{t^{1+1/r}}dt\leq 1+r$$
and 
$$\sum_{m=1}^{h-1} \frac{1}{m^{1-1/r}}\leq rh^{1/r}.$$
Thus by \eqref{63},
\begin{align*}
\sum_{m=1}^{h-1}P(F,m)r^{-1}
m^{-1-1/r} \leq & 2r^{-1}g(2^r h)(1+\log^*h^{1/r})\sum_{m=1}^{h-1}\left(3+\frac{e^{q(r)}m^{2/r}}{|D(F)|^{1/r(r-1)}}\right).m^{-1-1/r}\\
\leq & 6r^{-1}g(2^r h)(1+\log^*h^{1/r})\sum_{m=1}^{h-1}m^{-1-1/r}\\
+&2r^{-1}g(2^r h)(1+\log^*h^{1/r})\frac{e^{q(r)}}{|D(F)|^{1/r(r-1)}}\sum_{m=1}^{h-1}m^{-1+1/r}\\
\leq & g(2^r h)(1+\log^*h^{1/r})\bigg(6+6/r+\frac{2e^{q(r)}h^{1/r}}{|D(F)|^{1/r(r-1)}}\bigg)
\end{align*}
So by Lemma \ref{PFh}(ii), we have
$$
N(F,h)\leq g(2^r h)(1+\log^*h^{1/r})\left(6+\frac{2e^{q(r)}h^{2/r}}{|D(F)|^{1/r(r-1)}}+h^{1/r}\bigg(6+6/r+\frac{2e^{q(r)}h^{1/r}}{|D(F)|^{1/r(r-1)}}\bigg)\right)\\
$$
which proves the bound for $N(F,h).$
\end{proof}
\noindent
{\bf Remark.} Our main task is therefore to find the function $g(m)$ satisfying \eqref{Cgm}. Suppose we have 
\be\label{Gm}
P'(F,m) \leq G(m)
\ee
where $G(m)$ is an increasing function of $m.$ The above proof reveals that 
$$P(F,h)\leq 2G(2^rh) (1+\log^*h^{1/r})$$
and by Lemma \ref{PFh}(i),
\be\label{GNFh}
N(F,h)\leq 2G(2^rh)(1+\log^*h^{1/r})h^{1/r}.
\ee

\section{Height}
The {\it na\"ive} height of $F$ denoted by $H(F)$ is defined as
$$H(F)=\max(|a_0|,|a_1|,\cdots,|a_r|).$$
Two forms $F$ and $G$ are said to be {\it equivalent} if 
$$G=F_A\ \textrm{ for some} \ A\in SL_2(\mathbb Z).$$
Note that equivalent forms have the same discriminant. Further if $F$ and $G$ are equivalent, then 
$$P'(F,h)=P'(G,h).$$
Thus in order to estimate $P'(F,h),$ it is enough to estimate the corresponding quantity for any equivalent form. 
Suppose $(x,y)=(1,0)$ is a solution of \eqref{secondinequality}. Then
we call the form $F$ {\it normalized}. Note that $F(1,0)=a_0$, the leading coefficient of $F.$ 
If $(x,y)$ is a primitive solution of \eqref{secondinequality}, then there is an $A\in SL_2(\mathbb Z)$ such that $A^{-1}(x,y)=(1,0).$ Then
$$\frac{h}{2^r}\leq |F_A(1,0)|<h.$$
Thus if $F$ has a primitive solution to \eqref{secondinequality}, we can find an equivalent, normalized form. From now on we shall assume that
{\it $F$ is a normalized form.}

The Mahler measure $M(F)$ is defined as 
$$M(F)=|a_0|\prod_{i=1}^r\max(1,|\alpha_i|).$$
Mahler \cite{Mah1} showed that
\be\label{Ma1}
|D(F)|\leq r^rM(F)^{2r-2}.
\ee
Further from \cite{Mah2}
\be\label{Ma2}
\bigg(\binom{r}{[r/2]}\bigg)^{-1} H(F)\leq M(F)\leq (r+1)^{1/2}H(F).
\ee
We will say that a form $F$ is {\it reduced} if it is normalized and has smallest Mahler height among all normalized forms equivalent to $F.$

\section{Diophantine Approximation}
We have the following approximation of roots of $F(x,1)$ by $x/y$ for integers $x$ and 
$y.$ This result is a refinement, due to Stewart \cite[Lemma 3]{Stew} of an estimate of  Lewis and Mahler \cite{Lew-Mah}.
\begin{lemma}\label{Lew}
For any pair of coprime integers $(x,y)$ with $y\neq 0$ we have that
$$\min_{\alpha}\bigg|\frac{x}{y}-\alpha\bigg|\leq \frac{2^{r-1}r^{(r-1)/2}M(F)^{r-2}|F(x,y)|}{|D(F)|^{1/2}|y|^r}$$
where $\alpha$ runs through all the roots of the polynomial $f(x)=F(x,1).$
\end{lemma}
As a consequence of the above lemma, we get the following corollary.
\begin{corollary}\label{application}
Let $F$ be an irreducible form of degree $r$ with
$$|D(F)|\geq 2^{2r}r^r h^{2(r-1)}.$$ 
Suppose $|y|\geq (M(F)/h)^2.$ Then
$$|x|\leq 2M(F)|y|.$$
\end{corollary}
\begin{proof}
By \eqref{Ma1} and the given condition on $|D(F)|$ we see that $M(F)/h>1.$ Further Lemma \ref{Lew} implies that 
$$
\min_{\alpha}\bigg|\frac{x}{y}-\alpha\bigg|\leq \frac{M(F)^{r-2}}{2r^{1/2}h^{r-2}|y|^r}\leq \frac{1}{(M(F)/h)^{r+2}}<1.
$$
Let the minimum on the left hand side occur at $\alpha=\alpha_0.$ Then
$$\bigg|\frac{x}{y}\bigg|\leq |\alpha_0|+\frac{1}{2(M(F)/h)^{r+2}}\leq 2M(F).$$
This gives the assertion of the lemma.
\end{proof}
We recall here yet another approximation result due to Schmidt \cite[Lemmas 5 and  6]{Sch}.
\begin{lemma}\label{psi}
Suppose $F$ is a reduced and normalized form with $M=M(F)/h>e^{4r}.$
Let $(x,y)$ be a primitive pair with $|y|>0$ satisfying \eqref{secondinequality}. 
Then
there exist numbers $\psi_j=\psi_j(x,y), 1\leq j\leq r$ such that 
\be\label{2r}
\psi_j=0\ or\ \frac{1}{2r}\leq \psi_j\leq 1
\ee
and
$$\sum_{j=1}^r\psi_j\geq \frac{1}{2}$$
such that
\be\label{schmidtinequality}
\bigg|\alpha_i-\frac{x}{y}\bigg|<\frac{1}{(M(F)/h)^{\psi_i/2}y^2}\ whenever \ \psi_i>0,
1\leq i\leq r.
\ee
\end{lemma}
\noindent
{\bf Note.} (i) By the above lemma, it follows that to every primitive solution $(x,y)$ of  \eqref{secondinequality} with $y\neq 0,$ there exists an $r$-tuple $(\psi_1,\ldots,\psi_r),$ with at least some $\psi_j\neq 0$ satisfying \eqref{schmidtinequality}.  

\noindent
(ii) In fact, Schmidt proves the above lemma with the condition $M>100^r.$ This condition is not used in \cite[Lemma 5]{Sch} and used in \cite[Lemma 6]{Sch} where \eqref{schmidtinequality} is obtained. It is used to secure that $M^{1/2r}>7$ so that $M^{\psi_i}-\frac{7}{2}\geq \frac{1}{2} M^{\psi_i}\geq M^{\psi_i/2}.$ It is easy to see that this inequality holds even under the assumption that $M>e^{4r}.$

For any primitive solution $(x,y)$ to \eqref{secondinequality}, $x/y$  is
close to one of the roots $\alpha_i, 1\leq i\leq r$ of $f(x)$ which are close to the real line.
Schmidt \cite{Sch} showed in an ingeneous way, how these roots are not far from each other. He termed this phenomenon as {\it clustering of roots.} Below we gather Lemmas 9,10 and 11 from \cite{Sch} which describe this phenomenon.
\begin{lemma}\label{8r}
Suppose $f(x)\in \mathbb R[x]$ is a polynomial of degree $r\geq 3.$ Let $w$  roots of $f$ be in the square $[a-\epsilon,a]\times (0,\epsilon].$
Suppose $f(x)f'(x)<0$ for every real $x$ in the interval $[a+\epsilon,a+8r\epsilon].$ Then there are at least 
$$w/120(\log ew)$$
roots in the square $(a,a+8r\epsilon)\times(0,8r\epsilon).$ 

\end{lemma} 
This is \cite[Lemma 9]{Sch}. It is applied to get
\begin{lemma}\label{435}
Suppose $f(x)\in \mathbb R[x]$ is a polynomial of degree $r\geq 3.$ Let $\sigma=A+iB$
be a root of $f$ with $B>0$ and suppose that $f(x)f'(x)<0$ for real $x$ in the interior of the interval $[A, A+(9r)^tB)$ where $t$ is a positive integer. Then $f$ has at least 
$e^{\sqrt{t/68}}$ roots in the square
$$ [A,A+(9r)^tB)\times (0, (9r)^t B].$$

\end{lemma}
\begin{proof}
This lemma is just \cite[Lemma 10]{Sch} with  $e^{\sqrt{t/68}}$ roots instead of 
$e^{\sqrt{t}/16}=e^{\sqrt{t/256}}$ roots. We indicate here the beautiful argument of Schmidt. For any non-negative integer $m\leq t,$ let
$$
P_m=[A,A+(9r)^mB]\times (0,(9r)^mB].
$$
Suppose $w_m$ is the number of roots in $P_m.$ Then $w_0\geq 1.$ Let $m<t$ and
$a= A+(9r)^mB,w=w_m$ and $\epsilon=(9r)^m B.$ Then $P_m=[a-\epsilon,a]\times (0,\epsilon].$ Also 
$$a+8r \epsilon=A+(9r)^mB(1+8r)< A+(9r)^{m+1}B\leq A+(9r)^tB.$$ 
Hence 
$$(a+\epsilon, a+8r\epsilon]\subseteq [A,A+(9r)^tB)$$
and so $ff'<0$ on this interval. 
Apply Lemma \ref{8r} to conclude that $f$ has at least $w_m/120\log(ew_m)$ roots in the square
$$\bigg(A+(9r)^mB,A+(9r)^mB+8r(9r)^mB\bigg)\times \bigg(0, 8r(9r)^mB\bigg).$$
The above square is disjoint from $P_m$ and is contained in $P_{m+1}.$ Thus 
$$w_{m+1}\geq w_m+\frac{w_m}{120\log ew_m}.$$
So by induction on $m,$ we have $w_m\geq m+1, m\geq 1.$ This is improved by the following argument. 
Note that
$$ m+1\geq e^{\sqrt{m/68}}\ {\rm for}\ m\leq 4912.$$
So
$$w_m\geq e^{\sqrt{m/68}}\ {\rm for}\ m\leq 4912.$$
Let $m>4912.$ Now using
$$\log\bigg(1+\frac{1}{x}\bigg)\geq \frac{1}{x}\bigg(1-\frac{1}{2x}\bigg)\ {\rm for}\ x>1,$$
we get
$$
\log w_{m+1}\geq \log w_m+\log\bigg(1+\frac{1}{120\log ew_m}\bigg)>\log w_m+\frac{.99}{120\log ew_m}
$$ 
and so
\begin{align*}
(\log w_{m+1})^2\geq & (\log w_m)^2+\frac{1.98\log w_m}{ 120\log ew_m}\\
\geq& (\log w_m)^2+\frac{1.98 \log(m+1)}{120\log e(m+1)}\\
\geq & (\log w_m)^2+(0.0165)\left(1-\frac{1}{1+\log(m+1)}\right)\\
\geq & (\log w_m)^2+(0.0165)(0.8947)\\
\geq & (\log w_m)^2+0.01476\\
\geq & (0.01476)(m+1)
\end{align*}
giving
$$w_m\geq e^{(.1214)\sqrt{m}}>e^{\sqrt{m/68}}\ \textrm{for all} \ m>4912$$
and hence 
$$w_m\geq e^{\sqrt{m/68}}\ {\rm for\ all}\ m\geq 1.$$
Thus $P_t$ contains at least $e^{\sqrt{t/68}}$ roots. This proves the lemma.
\end{proof}
In \cite{Mue-Sch}, the above clustering property was applied to extract a set of roots of $F(x,1)$ having cardinality $\ll s$ and for which Lemma \ref{Lew} can be used.
See \cite[Lemmas 6-7]{Mue-Sch}.
We shall state and prove with a different value of $R,$ the results from these lemmas.
Let
{\color{black}\be\label{R}
R=e^{z(r)+5}
\ee}
where
$$z(r)=68\log^3 r+150\log^2 r+\log r$$
as given in \eqref{zr}.
\begin{lemma}\label{setS}
Let $f(z)\in \mathbb R[z]$ be a polynomial of degree $r.$ Suppose that $f(x)f'(x)\neq 0$  
for real $x\in I$ where $I$ is an interval $X_1<x<X_2$, or a half line $x<X_2$ or $x>X_1.$ Suppose there are $m\geq 1$ roots $\alpha_j=x_j+iy_j (1\leq j\leq m)$ with real parts $x_j\in I$ and $y_j\neq 0.$ Then there is a root $\alpha_\ell$ among these roots such that for every real $\xi\neq \alpha_i, 1\leq i\leq m$
\be\label{alphal11}
|\xi-\alpha_\ell|<R\min_{1\leq i\leq m}|\xi-\alpha_i|.
\ee

\end{lemma}
\begin{proof}
We may suppose that
$$f(x)f'(x)<0 \ \textrm{for}\ x\in I$$
by taking $f(-x)$ instead of $f(x),$ if necessary. Then $|f(x)|^2$ decreases in $I.$ So $I$ cannot contain arbitrarily large $x.$ Thus $I$ is either an interval of the type $X_1<x<X_2$ or $x<X_2.$ We first show that every root $\alpha_j=x_j+iy_j$ with $x_j\in I$ satisfies
\be\label{lbfory}
|y_j|\geq \frac{e^{2}}{R}(X_2-x_j), 1\leq j\leq m.
\ee
Suppose not. Then there is a root $z=A+iB$ with $A<X_2,B>0$ and
$$B< \frac{e^{2}}{R}(X_2-A).$$
Let 
$$t=\lceil 68\log^2 r\rceil.$$
Then
{\color{black}\begin{align*}
A+(9r)^tB<&A+e^{(68\log^2 r+1)(\log r+\log 9)}B\\
<&A+e^{z(r)+3}B\\
<&A+\frac{RB}{e^2}<X_2.
\end{align*}}
Hence by Lemma \ref{435}, $f$ has at least $e^{\sqrt{t/68}}>r$ roots which is a contradiction. Thus \eqref{lbfory} is satisfied.

Let $\alpha_\ell$ be a root for which $|y_\ell|$ is minimal among the $m$ roots $\alpha_1,\ldots, \alpha_m.$ Then
\begin{align*}
|\alpha_\ell-\alpha_i|\leq & |x_\ell-x_i|+|y_\ell-y_i|\\
\leq & |X_2-x_i|+|X_2-x_\ell|+2|y_i|\\
\leq & \frac{R}{e^{2}}|y_i|+\frac{R}{e^{2}}|y_\ell|+2|y_i|\\
\leq & 2\bigg(\frac{R}{e^{2}}+1\bigg)\ |y_i|.
\end{align*}
Thus for any real $\xi\neq \alpha_i, 1\leq i\leq m,$ we have
\begin{align}
|\xi-\alpha_\ell|\leq & |\xi-\alpha_i|+|\alpha_i-\alpha_\ell| \nonumber\\
\leq & |\xi-\alpha_i|+2\bigg(\frac{R}{e^{2}}+1\bigg)\ |y_i|\nonumber\\
\leq & \bigg(\frac{2R}{e^{2}}+3\bigg)\ |\xi-\alpha_i| \nonumber\\
\leq & R|\xi-\alpha_i|\nonumber
\end{align}
which proves \eqref{alphal11}.
\end{proof}

Define
$$
T(F)=
\begin{cases}
6s\ \textrm{if}\ F\ \textrm{has at\ most}\ s+1\ \textrm{non-zero coefficients}\\
3t+1\ \textrm{if} \ F\in C(t).
\end{cases}
$$
Suppose $F$ has no more than $s+1$ non-zero coefficients. Then $F_{A_j}\in C(4s-2).$ So $T(F_{A_j})=12s-5.$
As a result of Lemmas \ref{setS} and \ref{Lew} we get
\begin{corollary}\label{Lew1}
Let $(x,y)$ satisfy \eqref{secondinequality} with $y\neq 0.$ 
Then there exists a set $S$ of roots of $F(x,1)$ such that
\be\label{alphal1}
\bigg|\alpha_\ell-\frac{x}{y}\bigg|< R\min_{1\leq i\leq r}\bigg|\alpha_i-\frac{x}{y}\bigg|.
\ee
and
\be\label{alphal}
\bigg|\alpha_\ell-\frac{x}{y}\bigg|\leq \frac{R2^{r-1}r^{(r-1)/2}M(F)^{r-2}h}{|D(F)|^{1/2}|y|^r}
\ee
for some $\alpha_\ell\in S$ and $|S|\leq T(F).$
Similar result holds if $x\neq 0$ by replacing $x/y, \alpha_\ell, \alpha_i$ and $S$ with $y/x,\alpha_\ell^{-1}, \alpha_i^{-1}$ and $S^*$ where $\alpha^{-1}$ are the roots of $F(1,y)$ and $S^*$ is a set of roots of $F(1,y)$ with $|S^*|\leq T(F).$  
\end{corollary}

\begin{proof}
Let $f(z)=F(z,1).$
\vskip 1mm
\noindent
(i) Suppose $F(x,y)$  has at most $s+1$ non-zero coefficients, then $ff'$ has at most $4s-1$ real zeros (see Corollary \ref{realroots}). Hence the real numbers $x$ for which $f(x)f'(x)\neq 0$ fall into at most $4s$ intervals or half lines $I.$ Form a set $S$ consisting of roots $\alpha$ of $f$ as follows.
Either $\alpha$ is  real or $\alpha=\alpha_\ell=x_\ell+iy_\ell$ with $x_\ell\in I$ for some $I$ and $\alpha_\ell$ satisfies \eqref{alphal11}. The set $S$ satisfies
$$|S|\leq 2s+4s=6s.$$
The inequality \eqref{alphal1} is immediate from \eqref{alphal11}.

Now we apply Lemma \ref{Lew} to get \eqref{alphal}.
Let $F\in C(t).$ By the remarks preceding the Corollary \ref{realroots}, we find that $f$ has at most $t+1$ real roots and $ff'$ has at most $2t+1$ real roots. Arguing as in (i), we get
a set $S$ of roots $\alpha_\ell$ of $f(z)$  with $|S|\leq 3t+1$ such that  for every real $\xi,$ the inequality \eqref{alphal1} is valid and by Lemma \ref{Lew}, we get \eqref{alphal}.
\end{proof}
In many proofs below, we shall use the set $S.$ Similar proofs are valid with the set $S^*$ as well.
Schmidt defined the following quantity in order to separate the roots which cluster near the 
real line.
Let
$$M(F)>1.$$
For any $\alpha\neq 0,$ put
$$
\Phi(\alpha)=
\begin{cases}
0\ \textrm{if} \ |\textrm{Im} \ \alpha|\geq 1\\
\frac{\log|{\textrm Im}\ \alpha|^{-1} }{\log M(F)}\ \textrm{if}\ 0<|\textrm{Im}\ 
\alpha|< 1\\
\infty\ \textrm{if}\ |\textrm{Im}\ \alpha|=0
\end{cases}
$$
and
\be\label{defn_SF}
S(F)=1+\sum_{i=1}^r\min\bigg(1,\Phi(\alpha_i)\bigg).
\ee
Thus $S(F)$ counts all the real roots with weight 1 and the complex roots {\it close} to the real line with weight $\frac{\log|\textrm{Im}(\alpha_i)|^{-1} }{\log M(F)}.$ 
The roots which are not close to the real line are not counted in $S(F).$
Suppose
\be\label{SF1}
\min_{\alpha}\bigg|\alpha-\frac{x}{y}\bigg|=\bigg|\alpha_\ell-\frac{x}{y}\bigg|<1.
\ee
Then
$$M(F)^{-\Phi(\alpha_\ell)}=|\textrm{Im}(\alpha_\ell)|<1$$
giving $\Phi(\alpha_\ell)>1$ since $M(F)>1.$ The number of such $\alpha_\ell$ is at most $S(F).$

\section{ Thue-Siegel principle}\label{TSP}
The Thue-Siegel principle relates two {\it very good} approximations to an algebraic number.
We follow the notation of \cite[p.\, 74]{Bomb-Sch} for stating this principle. Let $t,\tau$ be positive numbers such that
$$ t<\sqrt{2/r}, \ \sqrt{2-rt^2}<\tau<t.$$
Put
$$\lambda=\frac{2}{t-\tau}\ \textrm{ and }\ A_1=\frac{t^2}{2-rt^2}\bigg(\log M(F)+\frac{r}{2}\bigg).$$
We say that a rational number $x/y$ is a very good approximation to an algebraic number $\alpha$ of degree $r$ if
$$\bigg|\alpha-\frac{x}{y}\bigg|<\bigg(4e^{A_1}H(x,y)\bigg)^{-\lambda}.$$ 
The following lemma summarizes the Thue-Siegel principle.
\begin{lemma}\label{Th-Sie}
If $\alpha$ is a root of $F(x,1)$ and $x/y,x'/y'$ are two very good approximations to $\alpha,$
then
$$\log (4e^{A_1})+\log H(x',y')\leq \delta^{-1}\big(\log(4e^{A_1})+\log H(x,y)\big)$$
where
$$\delta=\frac{rt^2+\tau^2-2}{r-1}.$$ 
\end{lemma}

For application, we choose

$$t=\sqrt{2/(r+a^2)},\tau=bt$$
with $0<a<b<1.$ Then 
\be\label{lambda}
\lambda=\frac{2}{(1-b)t}=\frac{\sqrt{2(r+a^2)}}{1-b}>\frac{\sqrt{2r}}{1-b},
\ee
\be\label{del}
A_1=\frac{1}{a^2}\bigg(\log M(F)+\frac{r}{2}\bigg)\ \textrm{ and }\ \delta=\frac{2(b^2-a^2)}{(r-1)(r+a^2)}.
\ee

\section{Counting {\it large} solutions of \eqref{secondinequality}}
We apply Thue--Siegel principle to  count large solutions. Here one can use the set $S$ constructed in Corollary \ref{Lew1} or the quantity $S(F)$  defined in \eqref{defn_SF} to derive the following proposition. As mentioned earlier, we take solutions $(x,y)$ with $y>0.$ 

\begin{lemma}\label{largesolutions}
Let
$$M=M(F)/h.$$
Assume that 
$$|D(F)|\geq R^22^{2r} r^{r-1}h^{2(r-1)}$$
and 
$\lambda<r.$
Let
$$
\nu=\frac{\bigg(3\log 2+r/2a^2+(1+1/a^2)\log h\bigg)(r-1)}{r\log 2-\frac{1}{2}\log r }
$$
and
$$
 \eta=\frac{(2+\nu+1/a^2)\lambda-2}{r-\lambda}.
$$
Then the number of primitive solutions $(x,y)$ of \eqref{secondinequality} with  
$y\geq M^2$ is at most
$${\min\bigg(T(F),S(F)\bigg)}\bigg(2+\bigg\lfloor\frac{\log \eta}{\log(r-1)}\bigg\rfloor+\bigg\lfloor\frac{1}{\log (r-1)}\log\bigg(\frac{(2+\nu+1/a^2)r-2}{\delta((2+\nu+1/a^2)\lambda-2)}\bigg) \bigg\rfloor\bigg).$$

\end{lemma}
\begin{proof}
Note that 
\be\label{MF}
M=M(F)/h\geq \bigg(\frac{2^r}{\sqrt{r}}\bigg)^{1/(r-1)}>1.
\ee
Let $T$ be the set of primitive solutions $(x,y)$ satisfying \eqref{secondinequality} with
$$y\geq M^2.$$ Let $S$ be the set given by Corollary \ref{Lew1}.
For $\alpha_{\ell}\in S$, let $T^{(\ell)}$ denote the set of $(x,y)\in T$ such that
\eqref{alphal} is satisfied.
Then $T=\cup T^{(\ell)}$. We shall estimate $T^{(\ell)}.$ 
From the assumption on $|D(F)|$ and \eqref{alphal}, we see that
$$\bigg|\alpha_\ell-\frac{x}{y}\bigg|\leq \frac{M^{r-2}}{2y^r}\leq 
\frac{1}{2M^{r+2}}<1.$$
Further $M(F)>1$ by \eqref{MF}. Hence as seen after \eqref{SF1}, we have
$$\Phi(\alpha_\ell)>1.$$
Arrange the solutions in $T^{(\ell)}$ as $(x_1,y_1),(x_2,y_2),\ldots$ such that $M^2\leq y_1\leq y_2\leq \ldots.$ Put $y_j=M^{1+\delta_j}$
so that $\delta_j\geq 1$ for $j\geq 1$ and
\be\label{goodappro}
\bigg|\alpha_\ell-\frac{x_j}{y_j}\bigg|\leq \frac{1}{2M^{2+r\delta_j}}.
\ee
Then
$$\frac{1}{y_jy_{j+1}}\leq \bigg|\frac{x_j}{y_j}-\frac{x_{j+1}}{y_{j+1}}\bigg| 
 \leq \bigg|\frac{x_j}{y_j}-\alpha_\ell\bigg|+\bigg|\frac{x_{j+1}}{y_{j+1}}-\alpha_{\ell}\bigg| <\frac{M^{r-2}}{y_j^r} $$
giving 
$$y_j^{r-1}\leq M^{r-2}y_{j+1}$$
i.e.,
$$M^{(1+\delta_j)(r-1)}\leq M^{r-1+\delta_{j+1}}.$$
Hence
$$\delta_{j+1}\geq (r-1)\delta_j.$$
It follows by induction that
\be\label{deltaj}
\delta_j \geq (r-1)^{j-1}\delta_1\geq (r-1)^{j-1}\ \textrm{for}\ j\geq 1.
\ee

Let us now see when $x_j/y_j$ will be a very good approximation to $\alpha_\ell$, as defined in Section \ref{TSP}. 
By \eqref{MF},
$$
M^\nu \geq e^{3\log 2+r/2a^2+(1+1/a^2)\log h}.
$$
Using Corollary \ref{application}, we first compute
\begin{align}
4e^{A_1}H(x_j,y_j)\leq & 8e^{r/2a^2}M(F)^{1+1/a^2}|y_j|\nonumber\\
=& 8e^{r/2a^2}M(F)^{2+\delta_j+1/a^2}h^{-1-\delta_j}\nonumber\\
=& 8e^{r/2a^2}M^{2+\delta_j+1/a^2}h^{1+1/a^2}\leq M^{2+\nu+\delta_j+1/a^2}\label{nudelta}.
\end{align}
Thus it follows from \eqref{goodappro} that $x_j/y_j$ is a very good approximation to $\alpha_\ell$ if
\begin{equation*}
2+r\delta_j\geq (2+\nu+\delta_j+1/a^2)\lambda
\end{equation*}
i.e., if
$$
\delta_j\geq \frac{(2+\nu+1/a^2)\lambda-2}{r-\lambda}=:\eta.
$$
By \eqref{deltaj}, we have $\delta_j\geq \eta$ for $j>J$ 
where
$$J=\bigg\lceil\frac{\log\eta}{\log(r-1)}\bigg\rceil.$$
So we see that $x_j/y_j$ is a very good approximation to $\alpha_\ell$ for $j\geq J+1.$ If $|T^{(\ell)}|> J$, let $t=|T^{(\ell)}|- J$.
Let us take two rationals, say $x_{J+1}/y_{J+1}$ and $x_{J+t}/y_{J+t}$ which are very good approximations to $\alpha_\ell$. Then by the Thue--Siegel principle,
$$\log (4e^{A_1})+\log H(x_{J+t},y_{J+t})\leq 
\delta^{-1}\bigg(\log (4e^{A_1})+\log H(x_{J+1},y_{J+1})  \bigg).$$
By \eqref{nudelta}, this implies that
$$1+\delta_{J+t}\leq \delta^{-1}(2+\nu+\delta_{J+1}+1/a^2).$$
Thus
\begin{align*}
(r-1)^{t-1}\leq \frac{\delta_{J+t}}{\delta_{J+1}}\leq &\delta^{-1}\bigg(1+\frac{2+\nu+1/a^2}{\delta_{J+1}}\bigg)\\
\leq &\delta^{-1}\bigg(1+\frac{2+\nu+1/a^2}{\eta}\bigg)\\
\leq &\delta^{-1}\bigg(\frac{(2+\nu+1/a^2)r-2}{(2+\nu+1/a^2)\lambda-2}\bigg).
\end{align*}
Taking logarithms on either side of the above inequality, we get
$$t\leq 1+
\bigg[ \frac{1}{\log (r-1)}\log\bigg(\frac{(2+\nu+1/a^2)r-2}{\delta((2+\nu+1/a^2)\lambda-2)}\bigg) \bigg].$$
Thus
$$|T^{(\ell)}|\leq 2+\bigg\lfloor\frac{\log \eta}{\log(r-1)}\bigg\rfloor+\bigg\lfloor\frac{1}{\log (r-1)}\log\bigg(\frac{(2+\nu+1/a^2)r-2}{\delta((2+\nu+1/a^2)\lambda-2)}\bigg) \bigg\rfloor.$$
Since
 there are at most $T(F)$ such $\alpha_\ell$ by Corollary \ref{Lew1} or $S(F)$ such $\alpha_\ell$ since $\Phi(\alpha_\ell)>1,$ we get the assertion of the lemma.
\end{proof}

We apply Lemma \ref{largesolutions} with specific values of $a$ and $b$ to get the following explicit result for large solutions. We have taken the same values of $a$ and $b$ as specified in \cite[p.252]{Div-Sar}, since these values yielded `good' results in that context and the problem here is similar. 
\begin{lemma}\label{largesolutions1}
Assume that $|D(F)|\geq R^22^{2r} r^{r-1}h^{2(r-1)}.$ Then
the number of primitive solutions $(x,y)$ of \eqref{secondinequality} with
$y\geq M^2$ is at most
$$2\min \bigg(T(F), S(F)\bigg)\ \bigg(15+\frac{1}{\log(r-1)}\log\log^*h^{1/r}\bigg).$$
\end{lemma}
\begin{proof}
The choices of $a$ and $b$ are as follows.
$$(a,b)=
\begin{cases}
(.5,.54)\ \textrm{if} \ r\geq 24\\
(.4,.48)\ \textrm{if}\ 9\leq r\leq 23\\
(.3,.36)\ \textrm{if} \ r=6,7,8\\
(.2,.24)\ \textrm{if} \ r=4,5\\
(.1,.15)\ \textrm{if}\  r=3.
\end{cases}
$$
These values will be used in the proof without any mention. We estimate the quantities $\lambda, \delta, \nu,\eta$ from their definitions (see Section 6 and Lemma \ref{largesolutions} ) as follows. These estimates are valid for $r\geq 3$ unless mentioned otherwise.
\be\label{lamb}
\lambda>\frac{\sqrt{2r}}{1-b};\ \ \lambda\leq \frac{1.64\sqrt{r}}{1-b},\ \delta>\frac{2(b^2-a^2)}{r^2}.
\ee
Next
\begin{align*}
\nu&\leq  \frac{r\bigg((3\log 2)/r+1/2a^2+(1+1/a^2)\log ^*h^{1/r}\bigg)}{\bigg(\log 2-\frac{\log r}{2r}\bigg)}\\
&\leq 2r\ \bigg(\log 2+1/2a^2+(1+1/a^2)\log^*h^{1/r}\bigg)\\
&\leq r\ (1.4+1/a^2)+2r(1+1/a^2)\log^*h^{1/r}.
\end{align*}
Thus
\be\label{h*}
2+\nu+1/a^2\leq 2r(1+1/a^2)(1+\log^*h^{1/r})\leq 4r(1+1/a^2)\log^*h^{1/r}.
\ee
Also
$$\nu\geq \frac{r(r-1)/2a^2}{r\log 2-(\log r)/2}\geq \frac{ r}{2a^2}.$$
Hence
$$(2+\nu+1/a^2)\lambda-2\geq \frac{r\lambda}{2a^2}\geq \frac{r^{3/2}}{\sqrt{2}a^2(1-b)}.$$
Further
$$r-\lambda= r-\frac{\sqrt{2r}}{1-b}\sqrt{1+a^2/r}.$$
Using the specified values of $a$ and $b$ we see that
$$r-\lambda\geq \alpha r$$
where
$$\alpha=
\begin{cases}
.36\ \textrm{for} \ r\geq 24\\
.08\ \textrm{for}\ 9\leq r\leq 23\\
.09\ \textrm{for}\ r=6,7,8\\
.06\ \textrm{for}\ r=4,5\\
.03\ \textrm{for}\ r=3.
\end{cases}
$$
Thus from the definition of $\eta$, \eqref{lamb}, \eqref{h*} and $r\geq 3,$ we get
\begin{align*}
&\frac{\log \eta}{\log(r-1)} \leq \frac{\log(2+\nu+1/a^2)+\log \lambda-\log (r-\lambda)}{\log(r-1)}\\
& \leq \frac{\log (4 (1.64)r^{1/2})}{\log(r-1)}+\frac{1}{\log(r-1)}\log\bigg( \frac{1+1/a^2}{\alpha (1-b)}\bigg)+ \frac{1}{\log(r-1)}\log \log^*h^{1/r}\\
&\leq 3.52+\frac{1}{\log(r-1)}\log\bigg( \frac{1+1/a^2}{\alpha (1-b)}\bigg)+ \frac{1}{\log(r-1)}\log \log^*h^{1/r}.
\end{align*}
Using the choices of $(a,b)$ and the values for $\alpha,$ we get 
\be\label{eta1}
\frac{\log \eta}{\log(r-1)} \leq 16+ \frac{1}{\log(r-1)}\log \log^*h^{1/r}.
\ee
Finally,
\begin{align*}
&\frac{1}{\log(r-1)}\log\bigg(\frac{(2+\nu+1/a^2)r-2}{\delta((2+\nu+1/a^2)\lambda-2)}\bigg)  \\
&\leq\frac{1}{\log(r-1)}\bigg(\log {\delta^{-1}}+\log\bigg(4r^2(1+1/a^2)\log^*h^{1/r}\bigg)-\log\bigg(\frac{r^{5/2}}{\sqrt{2}a^2(1-b)}\bigg)\bigg)\\
&\leq\frac{1}{\log(r-1)}\bigg(\frac{3}{2}\log r+\frac{3}{2}\log 2+\log\frac{(a^2+1)(1-b)}{(b^2-a^2)}\bigg)+\frac{\log\log^*h^{1/r}}{\log(r-1)}\\
& \leq 12+\frac{1}{\log(r-1)}\log\log^*h^{1/r}.
\end{align*}
Now the result follows from \eqref{eta1} and Lemma \ref{largesolutions}.
\end{proof}

\section{Counting small solutions of \eqref{secondinequality} using Bombieri-Schmidt \cite{Bomb-Sch}}

\subsection{Estimation of small solutions in terms of $r$ and $s$}
In this section we want to count the number of primitive solutions $(x,y)$ of \eqref{secondinequality} with $1\leq y\leq Y.$  Let us call this set $Sm(Y)=Sm(Y,F).$
Assume that $Sm(Y)\neq \emptyset.$ We refer to \cite[p.247-248]{Sch} for the following facts. Let 
$$L_i(x,y)=x-\alpha_iy\ {\rm for} \ 1\leq i\leq r.$$
For any two ${\bf x}=(x,y)$ and ${\bf x_0}=(x_0,y_0),$ let
$$D({\bf x},{\bf x_0})=xy_0-x_0y.$$
Let us assume that $F$ is reduced. Let $(x_0,y_0)$ and $(x,y)$ be two 
linearly independent primitive solutions of \eqref{secondinequality}. Suppose ${\bf x'}=(x',y')$ is such that $D({\bf x'},{\bf x})=1$ so that ${\bf x'},{\bf x}$ is a basis of $\mathbb Z^2.$ Set 
$$G(v,w)=F(v{\bf x}+w{\bf x'}).$$
Then $G\sim F$ and $G$ is normalized and 
\be\label{GVW}
G(v,w)=F({\bf x})\prod_{i=1}^r (v-\beta_iw), \beta_i\in \mathbb C, 1\leq i\leq r.
\ee
Then
from \cite[p.248-line 14]{Sch}, there exists an integer $m=m(x_0,y_0,x,y)$ such that
\be\label{Li}
\bigg|\frac{L_i(x_0,y_0)}{L_i(x,y)}\bigg|\geq \bigg(|m-\beta_i|-\frac{1}{2}\bigg)\ |D({\bf x},{\bf x_0})|-2.
\ee
Since $F$ is reduced, it is normalized and so $F(1,0)$ satisfies \eqref{secondinequality}. Taking $(x_0,y_0)=(1,0),$
we find
$y=|D({\bf x,x_0})|$ and from \eqref{Li}, we get
\be\label{Lii}
\bigg|\frac{1}{L_i(x,y)}\bigg|\geq \bigg(|m-\beta_i|-\frac{1}{2}\bigg)y-2
\ee
where $m=m(x,y).$
Define 
$$\chi_i=\big\{(x,y)\in Sm(Y): |x-\alpha_iy|\leq 1/(2y)\big\},\ 1\leq i\leq r.$$
Further let the elements of $\chi_i$ be given by
$${\bf x}^{(i)}_1,\ldots, {\bf x}^{(i)}_{\sigma_i}$$
ordered so that $y^{(i)}_1\leq \cdots\leq y^{(i)}_{\sigma_i}.$ Put ${\bf x}(i)={\bf x}^{(i)}_{\sigma_i}.$ Let
$$Z=\big\{{\bf x}(1), \ldots, {\bf x}(r)\big\}$$ 
and
$$\chi=Sm(Y)\backslash Z.$$
For the next lemma we follow \cite[Lemmas 4 and 5]{Bomb-Sch}.
\begin{lemma}\label{45}
Let $F$ be reduced.
Let $1\leq i\leq r$ and $(x,y)\in \chi.$ Then there exist $m=m(x,y)\in\mathbb Z$ and $\beta_i=\beta_i(x,y)\in \mathbb C$  satisfying the following properties.

(i) Suppose $(x,y)$ and $(x',y')$ are in $\chi_i$ with $y'\geq y.$ Then 
$$\frac{y'}{y} \geq \frac{2}{9}\max\big(1,|\beta_i-m|\big).$$

(ii) Suppose $(x,y)\not\in \chi_i.$ Then
$$|\beta_i-m|\leq 9/2.$$

(iii)$$\prod_{i=1}^r\max\big(1,|\beta_i-m|\big)\geq \frac{M(F)}{h}.$$
\end{lemma}
\begin{proof}
We have
\begin{align*}
1\leq |D({\bf x,x'})|&=|xy'-x'y|\\
&=\big|(x-\alpha_iy)y'-(x'-\alpha_iy')y\big|\\
&\leq y\big|L_i(x',y')\big|+y'\big|L_i(x,y)\big|\\
&\leq \frac{y}{2y'}+y'\big|L_i(x,y)\big|,
\end{align*}
so that 
\be\label{yi}
y'\big|L_i(x,y)\big|\geq 1/2.
\ee
Since $F$ is reduced, we get
by \eqref{Lii} and \eqref{yi} that
$$ y' \geq \frac{1}{2}\bigg(|m-\beta_i|-\frac{1}{2}\bigg)y-1$$
so that
$$\frac{y'}{y}\geq \max\bigg(1,\frac{1}{2}\big(|m-\beta_i|\big)-\frac{5}{4}\bigg).$$
It is easy to see that for any real $\xi\geq 0,$
$$\max\bigg(1,\frac{1}{2}\xi-\frac{5}{4}\bigg)\geq \frac{2}{9}\max(1,\xi).$$
Hence assertion $(i)$ follows.
\vskip 2mm
When $(x,y)\not\in \chi_i,$ we have 
$$|L_i(x,y)|>\frac{1}{2y}.$$
Using this in \eqref{Lii}, we get
$$|m-\beta_i|\leq 2+\frac{2}{y}+\frac{1}{2}\leq \frac{9}{2}$$
proving $(ii).$
\vskip 2mm
For $(iii),$ let
$$ G^*({\bf x})=G(x+my,y)=b_0\prod_{i=1}^r(x+(m-\beta_i)y)$$
where $G$ is given by \eqref{GVW}. Then $G^*$ is normalized and equivalent to $G$ hence to $F.$ Thus $|b_0|\leq h$ and since $F$ is reduced, we have $M( G^*)\geq M(F).$ But 
$$M(G^*)=|b_0|\prod_{i=1}^r \max(1,|\beta_i-m|)\leq  h\prod_{i=1}^r \max(1,|\beta_i-m|)$$
which proves $(iii).$
\end{proof}

As a consequence of Lemma \ref{45}, we get
\begin{lemma}\label{451}
Let $F$ be reduced with $M(F)>(9/2)^rh.$ Then
$$|\chi|\leq \frac{r\log Y}{\log M(F)-r\log(9/2)-\log h}.$$
\end{lemma}
\begin{proof}
Take any ${\bf x}^{(i)}_k\in \chi_i, 1\leq k<\sigma_i.$ Then
by Lemma \ref{45}(i),
$$\frac{2}{9}\max\bigg(1,\big|\beta_i({\bf x}^{(i)}_k)-m({\bf x}^{(i)}_k)\big|\bigg)\leq \frac{y^{(i)}_{k+1}}{y^{(i)}_k}.$$
Hence we get
$$\prod_{{\bf x}\in {\chi\cap \chi_i}}\frac{2}{9}\max\bigg(1,\big|\beta_i({\bf x})-m({\bf x})\big|\bigg)\leq 
\frac{y^{(i)}_{\sigma_i}}{y^{(i)}_1}\leq Y.$$
On the other hand, when ${\bf x}\not\in \chi_i,$ by Lemma \ref{45}(ii)
$$\frac{2}{9}\max\bigg(1,\big|\beta_i({\bf x})-m({\bf x})\big|\bigg)\leq 1.$$
Hence
$$\prod_{{\bf x}\in \chi}\frac{2}{9}\max\bigg(1,\big|\beta_i({\bf x})-m({\bf x})\big|\bigg)\leq Y.$$
This gives
$$\prod_{i=1}^r\prod_{{\bf x}\in \chi}\frac{2}{9}\max\bigg(1,\big|\beta_i({\bf x})-m({\bf x})\big|\bigg)\leq Y^r.$$
Note that the product on the left hand side is indeed
$$\prod_{{\bf x}\in \chi}\prod_{i=1}^r\frac{2}{9}\max\bigg(1,\big|\beta_i({\bf x})-m({\bf x})\big|\bigg).$$
Therefore, by Lemma \ref{45}(iii), we get
$$\bigg(\frac{M(F)}{(9/2)^rh}\bigg)^{|\chi|}\leq Y^r$$
which proves the result.
\end{proof}

Since $|Sm(Y)|\leq |\chi|+r$ it follows from Lemma \ref{451} that
\begin{corollary}\label{20s}
Let $F$ be reduced and $M(F)>(9/2)^rh.$ Then
$$|Sm(Y)|\leq r+\frac{r\log Y}{\log M(F)-r\log(9/2)-\log h}.$$
\end{corollary}
\begin{lemma}\label{smallsolutions}
Let $F$ be reduced.
Assume that 
$$|D(F)|\geq e^{10r(r-1)+2z(r)}h^{2(r-1)}.$$ Then
the number of primitive solutions $(x,y)$ of \eqref{secondinequality} with
$0<y\leq M^2$ is at most $4r.$
\end{lemma}
\begin{proof}
The assumption on $|D(F)|$ implies that
\begin{align*}
M(F)\geq& \frac{e^{5r+\frac{68\log^3r}{r-1}}h}{r^{\frac{r}{2(r-1)}}}\\
\geq& \frac{e^{5r}h}{\sqrt{r}}\geq e^{4.81r}h.
\end{align*}
Let $Y=M^2=(M(F)/h)^2.$
Note that
\begin{align*}
\log M(F)-r\log(9/2)-\log h&\geq \big(1-\log(9/2)/(4.81)\big)\big(\log M(F)-\log h\big)\\
&\geq (.687)\big(\log M(F)-\log h\big).
\end{align*}
Thus by the above corollary,
\begin{align*}
|Sm(Y)|&\leq r+\frac{2r\big(\log M(F)-\log h\big)}
{\log M(F)-r\log(9/2)-\log h}\\
&\leq r+\frac{2r\big(\log M(F)-\log h\big)}{(.687)\big(\log M(F)-\log h\big)}\\
&\leq r+(2.912)r\leq 4r.
\end{align*}
\end{proof}
{\color{black}In conclusion, by Lemmas \ref{largesolutions1}, \ref{smallsolutions} and including the solutions with $y<0$ and the solution $(1,0)$, we have the following proposition.
\begin{proposition}\label{redform}
Assume that $F$ is a reduced form of degree $r.$ Let 
$|D(F)|\geq e^{10r(r-1)+2z(r)}h^{2(r-1)}.$ 
Then the number of primitive solutions $(x,y)$ of \eqref{secondinequality} is at most
$$1+8r+4T(F)\bigg(15+\frac{1}{\log(r-1)}\log\log^*(2h^{1/r})\bigg).$$
\end{proposition}
\subsection{Final Assault for Theorem \ref{Main1}}
Let $F$ have at most $s+1$ non-zero coefficients.
By Lemma \ref{4s} and Corollary \ref{realroots} (iii), $F\in C(4s-2)$ and $F_{A_j}\in C_p(4s-2).$ 
In the choice of the prime $p$ in Section 3, we take 
$$q(r)=13+\frac{2z(r)}{r(r-1)}$$
so that by \eqref{lbp} we get
\begin{align*}
|D(F_{A_j})|&=p^{r(r-1)}|D(F)|\geq \frac{e^{q(r)r(r-1)}h^{2(r-1)}}{8^{r(r-1)}}\\
&\geq e^{(13-2.1)r(r-1)+2z(r)}h^{2(r-1)}\\
&\geq e^{10r(r-1)+2z(r)}h^{2(r-1)}.
\end{align*}
Among all matrices $B$ which are equivalent to $A_j$ and for which $F_B(x,y)$ is normalized, choose the one having minimal Mahler measure. Then $F_B$ is reduced and $|D(F_B)|$ also satisfies the above inequality. Further $T(F_B)=3(4s-2)+1<12s.$
So by Proposition \ref{redform} with $F=F_B,$ we get
$$P'(F_{A_j},h)=P'(F_{B},h)\leq 1+8r+48s\bigg(15+\frac{1}{\log(r-1)}\log\log^*h^{1/r}\bigg).$$
Take 
$$g(m)=1+8r+48s\bigg(15+\frac{1}{\log(r-1)}\log\log^*m^{1/r}\bigg).$$
Then
$$g(2^rh)\leq 1+8r+48s\bigg(15+\frac{1}{\log(r-1)}\log\log^*(2h^{1/r})\bigg)$$
and as $z(r)/r(r-1)\leq 45.4,$ we find
$$2e^{q(r)}\leq e^{105}.$$
Now Theorem \ref{Main1} follows from Lemma \ref{NCUg}. \qed

By Proposition \ref{redform} and \eqref{GNFh} with 
$$G(m)=1+8r+48s\bigg(15+\frac{1}{\log(r-1)}\log\log^*m^{1/r}\bigg).$$
we also get
\begin{corollary}
Suppose $F$ is reduced and having  at most $s+1$ non-zero coefficients with
$$|D(F)|\geq e^{10r(r-1)+2z(r)}h^{2(r-1)}.$$ 
Then 
$$N(F,h)\leq 2h^{1/r}(1+\log^*h^{1/r})\left(1+8r+24s\bigg(15+\frac{1}{\log(r-1)}\log\log^*(2h^{1/r})\bigg)\right).$$
\end{corollary}}
This may be compared with the results of  \cite{Sie} and \cite{Akh-Sara-Div} alluded to in the Introduction.
While estimating $|Sm(Y)|$ we saw earlier that the set
$$Z=\big\{{\bf x}(1), \ldots, {\bf x}(r)\big\}$$ 
is bounded trivially as $|Z|\leq r.$ This can be improved and it will be used in Section 10.
We follow the argument from \cite{Ben}. We show
{\color{black}\begin{lemma}\label{Z}
Let $F$ be reduced and $M(F)>(11/2)^rh. $
Let $S$ be the set given by Corollary \ref{Lew1}. Then
$$|Z|\leq 1+|S|+\frac{r\log(6R+5)-r\log 11}{\log M(F)-r\log(11/2)-\log h}$$
and
$$|Sm(Y)|\leq 1+|S|+\frac{r(\log Y+\log(6R+5)-\log 11)}{\log M(F)-r\log(11/2)-\log h}.$$
\end{lemma}}
\begin{proof}
Let us arrange the roots $\alpha_i$ of $F(x,1)$ as $\alpha_1,\alpha_2,\ldots,\alpha_t,\alpha_{t+1},\ldots,\alpha_r$ so that $\alpha_1,\alpha_2,\ldots, \alpha_t$ are all the roots in $S$ where $S$ is the set constructed in Corollary \ref{Lew1}. Thus $ |S|=t.$
Let us write 
$$Z=Z_1\cup Z_2$$
where $Z_1=\big\{{\bf x}(1), \ldots, {\bf x}(t)\big\}$ and $Z_2=\big\{{\bf x}(t+1), \ldots, {\bf x}(r)\big\}.$ Since we have already estimated $|\chi|$ in Lemma \ref{451} and we know $|Z_1|\leq |S|,$ it remains to estimate $|Z_2|$ and while doing so, we may assume that none of the ${\bf x}(i), t+1\leq i\leq r$ belongs to $\chi$ or $Z_1.$ Let us fix an ${\bf x}(i)\in Z_2.$ For simplicity, let us write ${\bf x}(i)=(u_i,v_i).$ Since $(u_i,v_i)\in \chi_i$ we have
\be\label{2vi}
|L_i(u_i,v_i)|\leq \frac{1}{2v_i}.
\ee
By the construction of $S,$ there exists $\alpha_\ell, 1\leq \ell\leq t $ such that
$$|L_\ell(u_i,v_i)|\leq R\min_{1\leq j\leq r}|L_j(u_i,v_i)|\leq R|L_i(u_i,v_i)|$$
giving
\be\label{2Rvi}
|L_i(u_i,v_i)|\geq \frac{1}{R}|L_\ell(u_i,v_i)|> \frac{1}{2Rv_i}
\ee
since $(u_i,v_i)\not \in \chi_\ell.$ We take $(x_0,y_0)\in Sm(Y)$ with the least $y$ value.
We assume that $(x_0,y_0)\neq (u_i,v_i).$
Then by \eqref{Li}, there exists $m=m(x_0,y_0,u_i,v_i)\in \mathbb Z$ and $\beta_i=\beta_i(u_i,v_i)\in \mathbb C$ such that
$$\left|\frac{L_i(x_0,y_0)}{L_i(u_i,v_i)}\right|\geq \left(|m-\beta_i|-\frac{1}{2}\right)\ |x_0v_i-y_0u_i|-2.$$
Using \eqref{2vi}, \eqref{2Rvi} and $|x_0v_i-y_0u_i|\geq 1$ we have
\begin{align}
|m(x_0,y_0,u_i,v_i)-\beta_i(u_i,v_i)|&\leq \frac{2v_iR|x_0-\alpha_iy_0|}{|x_0v_i-y_0u_i|}+\frac{5}{2}\nonumber\\
&\leq \frac{2R\big|\frac{x_0}{y_0}-\alpha_i\big|}{\big|\frac{x_0}{y_0}-\frac{u_i}{v_i}\big|}+\frac{5}{2}\nonumber\\
&\leq \frac{2R\bigg(\big|\frac{u_i}{v_i}-\alpha_i\big|+\big|\frac{x_0}{y_0}-\frac{u_i}{v_i}\big|\bigg)}{\big|\frac{x_0}{y_0}-\frac{u_i}{v_i}\big|}+\frac{5}{2}\nonumber\\
&\leq 2R\bigg(1+\frac{y_0}{2v_i}\bigg)+\frac{5}{2}\nonumber\\
&\leq (6R+5)/2.\label{R52}
\end{align}
Now let us take $(u_j,v_j)\in Z_2$ with $(u_j,v_j)\neq (u_i,v_i).$ Since $(u_j,v_j)\not \in \chi,$ by assumption, we have
$$L_i(u_j,v_j)>\frac{1}{2v_j}$$
and hence
\begin{align*}
|m(x_0,y_0,u_j,v_j)-\beta_i(u_j,v_j)|&\leq  \bigg(\frac{L_i(x_0,y_0)}{L_i(u_j,v_j)}+2\bigg)\frac{1}{|u_jy_0-v_jx_0|}+\frac{1}{2}\\
&\leq \frac{L_i(x_0,y_0)}{L_i(u_j,v_j)|u_jy_0-v_jx_0|}+\frac{5}{2}\\
&= \frac{\big|\frac{x_0}{y_0}-\alpha_i\big|}{v_j^2\big|\frac{u_j}{v_j}-\alpha_i\big|\big|\frac{u_j}{v_j}-\frac{x_0}{y_0}\big|}+\frac{5}{2}
\end{align*}
\begin{align*}
&\leq \frac{\big|\frac{u_j}{v_j}-\alpha_i\big|+\big|\frac{u_j}{v_j}-\frac{x_0}{y_0}\big|}{v_j^2\big|\frac{u_j}{v_j}-\alpha_i\big|\big|\frac{u_j}{v_j}-\frac{x_0}{y_0}\big|}+\frac{5}{2}\\
&\leq \frac{1}{v_j^2\big|\frac{u_j}{v_j}-\alpha_i\big|}+\frac{1}{v_j^2\big|\frac{u_j}{v_j}-\frac{x_0}{y_0}\big|}+\frac{5}{2}\\
&\leq \frac{2}{v_j}+\frac{y_0}{v_j}+\frac{5}{2}\leq \frac{11}{2}.
\end{align*}
In other words, for any $(u_j,v_j)\in Z_2$ and $\neq (u_i,v_i)$ we have
$$\frac{2}{11}|m(x_0,y_0,u_j,v_j)-\beta_i(u_j,v_j)|\leq 1.$$
Thus by \eqref{R52}, we get
$$\prod_{(u,v)\in Z_2}\frac{2}{11}\max(1,|m(x_0,y_0,u,v)-\beta_i(u,v)|)\leq \frac{6R+5}{11}$$
and so
$$\prod_{1\leq i\leq r}\prod_{(u,v)\in Z_2}\frac{2}{11}\max(1,|m(x_0,y_0,u,v)-\beta_i(u,v)|)\leq \big(\frac{6R+5}{11}\big)^r.$$
Interchanging the products on the left hand side and using Lemma  \ref{45}(iii) we get
$$\bigg(\big(\frac{2}{11}\big)^r \frac{M(F)}{h}\bigg)^{|Z_2|}\leq \big(\frac{6R+5}{11}\big)^r.$$
Thus
$$|Z_2|\leq \frac{r\log(6R+5)-r\log 11}{\log M(F)-r\log(11/2)-\log h}.$$
{\color{black}Now the assertion of the lemma follows from Lemma \ref{451} and including the solution $(x_0,y_0).$}
\end{proof}
\noindent
{\bf Remark. } The number of primitive solutions with $1\leq x\leq Y$ can also be estimated as in Lemma \ref{Z} and Corollary \ref{20s} by writing
$$F(x,y)=a_r(y-\gamma_1x)\cdots(y-\gamma_rx)$$
and taking
$$L_i(x,y)=y-\gamma_ix.$$
\section{Counting {\it small} solutions of \eqref{secondinequality} using  Schmidt \cite{Sch}}
\subsection{Estimation of small solutions in terms of $S(F)$}
\begin{lemma}\label{12SF}
Let $M>e^{4r}.$
The number of primitive solutions $(x,y)$ of \eqref{secondinequality} with
$0<y\leq M^2$ is at most
$20(S(F)-1).$
\end{lemma}
\begin{proof}
Let $S$ be the set of primitive solutions $(x,y)$ satisfying \eqref{secondinequality} with
$$0<y\leq M^2.$$
We estimate $|S|.$ 
As noted after Lemma \ref{psi}, corresponding to every $(x,y)\in S,$ there exists a tuple 
$(\psi_1,\ldots, \psi_r)$ with some $\psi_j=\psi_j(x,y)\neq 0$ satisfying \eqref{schmidtinequality}. Fix $1\leq j\leq r.$  Let
$$S^{(j)}=\{ (x,y)\in S: \psi_j(x,y)>0\}.$$
Then by \eqref{2r},
\be\label{S}
|S|=\sum_{(x,y)\in S} 1\leq 2\sum_{(x,y)\in S} \sum_{j=1}^r \psi_j(x,y)=
2\sum_{j=1}^r\sum_{(x,y)\in S^{(j)}}\psi_j(x,y).
\ee
Thus it is enough to estimate
$$\sum_{(x,y)\in S^{(j)}}\psi_j(x,y)$$
where we may assume that $S^{(j)}\neq \emptyset.$
Suppose $(x_1,y_1),\ldots, (x_\nu,y_\nu)$ are the elements of $S^{(j)}$ arranged such that $y_1\leq y_2\leq \cdots\leq y_\nu.$
By Lemma \ref{psi}, 
\be\label{middle}
\big|{\textrm{Im}}\ \alpha_j\big|\leq \bigg|\alpha_j-\frac{x_t}{y_t}\bigg|<\frac{1}{M^{\psi_j(x_t,y_t)/2}y_t^2}<1\ {\rm for}\ 1\leq t\leq \nu.
\ee
On the other hand, $|{\textrm{Im}}\ \alpha_j|\geq M^{-\Phi(\alpha_j)}.$ Thus
$$\psi_j(x_t,y_t)<2\Phi(\alpha_j)$$
giving
\be\label{psit}
\psi_j(x_t,y_t)\leq \min(1,2\Phi(\alpha_j))
\ee
 and
$$y_t^2<M^{\Phi(\alpha_j)-\psi_j(x_t,y_t)/2}<M^{\Phi(\alpha_j)}\ \textrm{for} \ 1\leq t\leq \nu.$$
Let $\xi_j=\xi(\alpha_j)=\min(2,\Phi(\alpha_j)/2).$ Then
\be\label{xi}
y_t\leq M^{\xi_j}\ \textrm{for}\ 1\leq t\leq \nu.
\ee
Let us now consider two distinct primitive solutions $(x_t,y_t)$ and $(x_{t+1},y_{t+1}).$ We know that
$$x_ty_{t+1}-x_{t+1}y_t\neq 0.$$  
We use the middle inequality in \eqref{middle} and $\psi_j\geq 1/2r$ to get
\begin{align*}
1\leq |x_ty_{t+1}-x_{t+1}y_t|=&|y_{t+1}(x_t-\alpha_j y_t)-y_t(x_{t+1}-\alpha_j y_{t+1})|\\
\leq & \frac{y_{t+1}}{y_t M^{\psi_j(x_t,y_t)/2}}+\frac{y_{t}}{y_{t+1} M^{\psi_j(x_{t+1},y_{t+1})/2}}\\
\leq & \frac{y_{t+1}}{y_t M^{\psi_j(x_t,y_t)/2}}+\frac{1}{M^{1/4r}}.
\end{align*}
Since $M>e^{4r},$ it follows that
$$y_{t+1}>(.632) M^{\psi_j(x_t,y_t)/2}y_t>M^{\psi_j(x_t,y_t)\big((.5)-(.46)(2r)/\log M\big)}y_t>M^{(.27)\psi_j(x_t,y_t)}y_t.$$
 Iterating the above inequality, we obtain
$$y_\nu\geq M^{(.27)\big(\psi_j(x_1,y_1)+\cdots+\psi_j(x_{\nu-1},y_{\nu-1})\big)}$$
which together with \eqref{xi} gives
$$\psi_j(x_1,y_1)+\cdots+\psi_j(x_{\nu-1},y_{\nu-1})\leq (3.71) \xi_j\leq (7.42)\min(1,\Phi(\alpha_j)).$$
Thus by \eqref{psit}, we get
\begin{align*}
\sum_{(x,y)\in S^{(j)}}\psi_j(x,y)=&\psi_j(x_1,y_1)+\cdots+\psi_j(x_{\nu-1},y_{\nu-1})+\psi_j(x_\nu,y_\nu)\\
\leq & (7.42) \min(1,\Phi(\alpha_j))+2\min(1,\Phi(\alpha_j))\\
<&10\min(1,\Phi(\alpha_j)).
\end{align*}
From \eqref{S}, it follows that
$$|S|\leq 20 \sum_{j=1}^r\min(1,\Phi(\alpha_j))$$
which proves the result.
\end{proof}
\subsection{Estimation of $S(F)$}
We first derive a simple consequence of Lemma \ref{435} which is \cite[Lemma 11]{Sch}.
\begin{lemma}\label{cor435}
Let $f\in \mathbb R[x]$ be a polynomial of degree $r\geq 3$ and $\epsilon>0$ and let $X_1<X_2.$ Suppose there are $u$ roots of $f$ in $[X_1,X_2]\times (0,\epsilon]$ with $ff'\neq 0$ in $(X_1,X_2).$ Then there are at least $u$ roots of $f$ in $[X_1,X_1+R\epsilon]\times (0,R\epsilon) .$
\end{lemma}
\begin{proof}
The assertion is obvious if $X_2-X_1\leq R\epsilon.$ Suppose that $X_2-X_1>R\epsilon.$ Let $t$ be the smallest integer with $e^{\sqrt{t/68}}\geq u.$ Then
$$68\log^2 u\leq t<68\log^2 u+1\leq 68\log^2 r+1$$
and
$$(9r)^t \epsilon<(9r)^{68\log^2 r+1}\epsilon<R\epsilon<X_2-X_1.$$
We may assume without loss of generality that $A=X_1$ and $B=\epsilon$ in Lemma \ref{435} and $ff'<0$ in $(X_1,X_2).$ Then there are at least $u$ roots in
$[X_1,X_1+(9r)^t\epsilon)\times (0,(9r)^t\epsilon]\subseteq [X_1,X_1+R\epsilon]\times (0,R\epsilon).$ This proves the lemma.
\end{proof}
We apply the above lemma to get the following result which can be found in \cite[Lemma 12, Corollary]{Sch}.

\begin{lemma}\label{qr}
Let $f\in \mathbb R[x]$ be a polynomial of degree $r\geq 3.$
Suppose that $ff'$ has $q-1$ real roots with $q\geq 1.$ Further suppose that  $M(f)>e^{4r}$ and 
 $$\frac{1+\log R}{r}\leq 1.$$
Choose a real number $\phi$ such that
$$\frac{1+\log R}{r}\leq \phi\leq 1$$
and put
$$f_\phi=\big\{\alpha: f(\alpha)=0 \ and\ 0<|{\rm Im}\ \alpha|\leq M(f)^{-\phi}\big\}.$$
Then
$$\big|f_\phi\big|\leq \bigg(\frac{6qr}{\phi}\bigg)^{1/2}.$$
\end{lemma} 
\begin{proof}
The assertion of the lemma is trivial if $q>r.$ So take $q\leq r.$  Let the $q-1$ real roots of $ff'$ be arranged as
$$X_1<X_2<\cdots<X_{q-1}.$$
This divides the real line into $q$ intervals
$$I_0=(-\infty, X_1],I_1=(X_1,X_2],\ldots,I_{q-2}= (X_{q-2},X_{q-1}],
I_{q-1}=(X_{q-1},\infty)$$
over each of the corresponding open interval, $ff'\neq 0.$
Divide the set $f_\phi$ into $q$ subsets 
$$f_j=\big\{\alpha\in f_\phi: \textrm{Re}(\alpha)\in I_j\big\}, 0\leq j\leq q-1 .$$
Note that $|f_j|$ is even.
Applying Lemma \ref{cor435} to each $I_j$ with $\epsilon=M(f)^{-\phi},$ we find that there are
at least $|f_j|/2$ roots of $f$ in the rectangle $[X_j, X_j+RM(f)^{-\phi}]\times (0,RM(f)^{-\phi}).$ Taking into account their complex conjugates as well, there are at least $|f_j|$ roots of $f$ in the rectangle $[X_j, X_j+RM(f)^{-\phi}]\times (-RM(f)^{-\phi},RM(f)^{-\phi})$ with mutual distance 
$$\leq \sqrt{5}RM(f)^{-\phi}<M(f)^{-3\phi/4}$$
since
$$M(f)^{\phi/4}>\sqrt{5}R$$
by the choice of $\phi$ and $M(f)>e^{4r}.$ Let us denote the set of these $|f_j|$ roots as $f_j'.$ Then $|f_j'|= |f_j|, 0\leq j\leq q-1.$
Consider
$$P=\prod_{j=0}^{q-1}\prod_{(\alpha,\beta)}\big|\alpha-\beta\big|$$
where the inside product is taken over pairs $(\alpha,\beta)$ with $\alpha$ and $\beta\in f_j'.$ Then
$$P>M(f)^{1-r}2^{-r(r-1)/2}$$
by \cite[Lemma 8, (7.3)]{Sch}. On the other hand, 
$$P\leq \prod_{j=0}^{q-1} M(f)^{-\binom{|f_j'|}{2}\frac{3\phi}{4}}=
M(f)^{-\frac{3\phi}{4}\sum_{j=0}^{q-1} \binom{|f_j'|}{2}}.$$
We may assume without loss of generality that all $|f_j'|>0$ hence $|f_j'|\geq 2, 0\leq j\leq q-1.$
Then using the inequality
$$x_1^2+\cdots+x_q^2\geq (x_1+\cdots+x_q)^2/q$$
for positive numbers $x_1,\ldots,x_q$ and
$\sum_{j=0}^{q-1}|f_j'|\geq 2q$ we get
\begin{align*}
\sum_{j=0}^{q-1}\binom{|f_j'|}{2}=&\frac{1}{2}\sum_{j=0}^{q-1}|f_j'|^2-\frac{1}{2}\sum_{j=0}^{q-1}|f_j'|\\
\geq & \frac{1}{2q}\bigg(\sum_{j=0}^{q-1}|f_j'|\bigg)^2\bigg(1-\frac{1}{2}\bigg)\\
=& \frac{1}{4q}\bigg(\sum_{j=0}^{q-1}|f_j'|\bigg)^2\\
=& \frac{1}{4q}\bigg(\sum_{j=0}^{q-1}|f_j|\bigg)^2=\frac{|f_\phi|^2}{4q}.
\end{align*}
Hence
$$2^{-r(r-1)/2}M(f)^{1-r}<P\leq \prod_{j=0}^{q-1} M(f)^{-\binom{|f_j'|}{2}{\frac{3\phi}{4}}}<M(f)^{-3\phi|f_\phi|^2/(16q)}$$
giving
$$|f_\phi|<\bigg(\frac{6qr}{\phi}\bigg)^{1/2}.$$
\end{proof}
\vskip 2mm
Now we come to the estimation of $S(F).$ Let $f(x)=F(x,1)$ and $\alpha_1,\ldots, \alpha_r$ be the roots of $f(x).$ Then $M(f)=M(F).$  Let us assume that
$$\frac{1+\log R}{r}\leq 1$$
and $\phi$ be such that
$$\frac{1+\log R}{r}\leq \phi\leq 1.$$
Divide the roots $\alpha_i$'s of $f(x)$ into three parts:
\begin{align*}
&(i)\ T_1=\big\{\alpha_i:\Phi(\alpha_i)\geq 1\big\},\\
 &(ii)\ T_2=\big\{\alpha_i:\Phi(\alpha_i)\leq \phi \big\},\\ 
&(iii)\ T_3=\big\{\alpha_i:\phi< \Phi(\alpha_i)<1\big\}.
\end{align*}
First we consider $T_1.$ For any $\alpha_i\in T_1$ we have $\big|\textrm{Im}\ \alpha_i\big|\leq M(F)^{-1}.$
Applying Lemma \ref{qr} with $\phi=1,$ we get
$$|T_1|\leq (6qr)^{1/2}.$$
Hence 
$$\sum_{\alpha_i\in T_1}\min\big(1,\Phi(\alpha_i)\big)=\sum_{\alpha_i\in T_1}1=|T_1|\leq (6qr)^{1/2}.$$
Let us fix $\phi=\frac{1+\log R}{r}$
and consider $T_2.$ 
Then
$$\sum_{\alpha_i\in T_2}\min\big(1,\Phi(\alpha_i)\big)\leq \sum_{\alpha_i\in T_2}\Phi(\alpha_i)\leq \phi r=1+\log R.$$
Lastly, take $T_3.$ Then $\Phi(\alpha_i)<1$ for any $\alpha_i\in T_3.$  Let us assume without loss of generality that $T_3=\big\{\alpha_1,\alpha_2,\ldots, \alpha_{|T_3|}\big\}.$
Arrange the $\alpha_i$'s in $T_3$ so that
$$\Phi(\alpha_1)\geq \Phi(\alpha_2)\geq \cdots \geq \Phi(\alpha_{|T_3|}).$$
Take any $\alpha_i.$
Then for any $j<i,$
$$\big|\textrm{Im}\ \alpha_j\big|=M^{-\Phi(\alpha_j)}\leq M^{-\Phi(\alpha_i)}.$$
Applying Lemma \ref{qr} with $\phi=\Phi(\alpha_i),$ we get
$$i\leq \big(6qr/\Phi(\alpha_i)\big)^{1/2}$$
giving
$$\Phi(\alpha_i)\leq 6qr/i^2.$$
Thus
$$\sum_{\alpha_i\in T_3}\min\big(1,\Phi(\alpha_i)\big)\leq \sum_{i=1}^\infty \min\bigg(1,\frac{6qr}{i^2}\bigg).$$
Now
\begin{align*}
\sum_{i=1}^\infty \min\bigg(1,\frac{6qr}{i^2}\bigg)\leq &(6qr)^{1/2}+
\sum_{i>(6qr)^{1/2}}^\infty \frac{6qr}{i^2}\\
\leq & (6qr)^{1/2}+\int_{\sqrt{6qr}-1}^\infty \frac{6qr}{x^2}dx=(2.5)(6qr)^{1/2}.
\end{align*}
So
$$\sum_{\alpha_i\in T_3}\min\{1,\Phi(\alpha_i)\}\leq (2.5)(6qr)^{1/2}.$$
From all the three contributions we get
$$S(F)-1=\sum_{\alpha_i}\min\big(1,\Phi(\alpha_i)\big)\leq 1+\log R+(3.5)(6qr)^{1/2}.$$
Using $q\leq 8s,$ we get
\begin{lemma}\label{SF}
Let $f(x)=F(x,1).$
Suppose $ff'$ has $q-1$ real roots with $q\geq 1.$ Let $M>e^{4r}$ and
$$\frac{1+\log R}{r}\leq 1.$$
Then
$$S(F)\leq 7+z(r)+25(sr)^{1/2}.$$
\end{lemma}
{\color{black}Combining Lemmas \ref{largesolutions1}, \ref{12SF},\, \ref{SF} and counting solutions with $y<0$ and $(1,0),$ we get the following proposition.
\begin{proposition}\label{105}
Assume that $F$ is a reduced form of degree $r$ with at most $s+1$ non-zero coefficients. Let 
$$|D(F)|\geq e^{10r(r-1)+2z(r)}h^{2(r-1)}$$
and
$$\frac{1+\log R}{r}\leq 1.$$  
Then the number of primitive solutions $(x,y)$ of \eqref{secondinequality} does not exceed
$$ 1+40\big(6+z(r)+25(sr)^{1/2}\big)+4T(F)\bigg(15+\frac{1}{\log(r-1)}\log\log^*h^{1/r}\bigg).$$

\end{proposition}
\noindent
{\bf Remark.} Note that  the condition on $r$ is satisfied for $r\geq 132434.$
\subsection{Final Assault for Theorem \ref{Main2}}
We follow the argument as in Section 8.2, use Proposition \ref{105} in place of Proposition \ref{redform} to get 
\begin{align*}
   & P'(F_{A_j},h)=P'(F_{B},h)\\
    &\leq 1+40\big(6+z(r)+25(sr)^{1/2}\big)+4T(F_B)\bigg(15+\frac{1}{\log(r-1)}\log\log^*h^{1/r}\bigg)
\end{align*}
where $T(F_B)<12s.$
We take 
$$g(m)=1+40\big(6+z(r)+25(sr)^{1/2}\big)+48s\bigg(15+\frac{1}{\log(r-1)}\log\log^*m^{1/r}\bigg).$$
Now Theorem \ref{Main2} follows from Lemma \ref{NCUg}. }
\qed



\section{Large, medium and small solutions \`{a} la Mueller and Schmidt}
We shall now proceed to prove Theorem 1.3 which is based on the method of 
Mueller and Schmidt \cite{Mue-Sch}.  We assume throughout that $r\geq 3s.$
We define three quantities as in \cite{Mue-Sch} by which the solutions of \eqref{secondinequality} are divided. We use $R$ as in \eqref{R}, $\lambda,\delta$ and $A_1$ as in Section 6.
Let
$$Y_L=\big(2^{r+1}r^{r/2}M(F)^rRh\big)^{1/(r-\lambda)} \big(4e^{A_1}\big)^{\lambda/(r-\lambda)}$$
and
$$Y_S=\big((e^6s)^rR^{2s}h\big)^{1/(r-2s)}.$$
Divide the primitive solutions $(x,y)$ of \eqref{secondinequality} into three sets according as
\begin{enumerate}[(i)]
\item $\max\big(|x|,|y|\big)>Y_L;$
\item $\max\big(|x|,|y|\big)\leq Y_L\ {\rm and}\ \min\big(|x|,|y|\big)>Y_S;$
\item $\min\big(|x|,|y|\big)\leq Y_S.$
\end{enumerate}
Denote the number of primitive solutions in these sets as $P_{lar}(Y_L)$, $P_{med}(Y_S,Y_L)$ and $P_{sma}(Y_S),$ respectively.

\subsection {Large solutions}

Let $(x,y)$ be a primitive solution of \eqref{secondinequality} with $\max(|x|,|y|)>Y_L.$
Put
\be\label{minalpha}
\min_{\alpha\in S}\left|\alpha-\frac{x}{y}\right|=\left|\alpha_\ell-\frac{x}{y}\right|
\ee
where $S$ is the set of roots of $F(x,1)$ given by Corollary \ref{Lew1}.
{\color{black}In \cite[Section III]{Bomb-Sch}, it was shown that the number of such primitive solutions with $\max(|x|,|y|)>Y_L$ and satisfying \eqref{minalpha} is 
$$< 2+\frac{\log(\delta^{-1}(\lambda-2)^{-1})}{\log(r-1)}$$
by taking $(x,y)$ and $(-x,-y)$ as one solution.
See also \cite[Lemma 9]{Mue-Sch}.
We use the estimates for $\lambda$ and $\delta$ from \eqref{lamb}
for $r\geq 24$ and their exact values from \eqref{lambda} and \eqref{del} for $r<24$ to  get the above value $<11.$ Thus in our case where we count $(x,y)$ and $(-x,-y)$ as two different solutions we find that the number of primitive solutions with $\max(|x|,|y|)>Y_L$ and satisfying \eqref{minalpha} is 
\be\label{lar}
P_{lar}(Y_L)\leq 22|S|\leq 132s
\ee
since $T(F)\leq 6s$ as $F$ has at most $s+1$ non-zero coefficients.}
\subsection { Small solutions}
We know $F\in C(4s-2).$ We argue as in Section 8.2, with the same $q(r)$ to find a matrix $B$ equivalent to $A_j$ with $F_B$ reduced. Further since
$$|D(F_B)|\geq e^{10r(r-1)+2z(r)}h^{2(r-1)}$$
we get
$$|M(F_B)|\geq e^{5r+\frac{z(r)}{r-1}}h.$$
Hence
$$\log |M(F_B)|-r\log(11/2)-\log h\geq 5r+\frac{z(r)}{r-1}-r\log(11/2)\geq 3.29r.$$
Also $|S|\leq T(F_B)\leq 12s-5.$
Apply Corollary \ref{Z} with $Y=Y_S, \log R=z(r)+5$ and $r\geq 3s$ to get
\begin{align}
Sm(Y_S,F_B)\leq &1+|S|+\frac{r\bigg(\log Y_S+\log(6R+5)\bigg)}{\log |M(F_B)|-r\log (11/2)-\log h}\nonumber\\
\leq & 12s-4+\frac{r\bigg(\frac{r}{r-2s}\big(6+\log s+\frac{2s}{r}\big(z(r)+5\big)+\log h^{1/r}\big)+z(r)+5+\log 7\bigg)}{3.29r}\nonumber\\
\leq & 12s-4+\frac{35+3\log s+3z(r)+3\log h^{1/r}}{3.29}\nonumber\\
\leq & 12s+z(r)+\log h^{1/r}.\nonumber
\end{align}
{\color{black}Thus the number of primitive solutions of \eqref{secondinequality} with $F=F_B$ and $0<y\leq Y_S$ is 
$$\leq 12s+z(r)+\log h^{1/r}.$$
Similar bound is valid when $0<x\leq Y_S.$
Note that in \eqref{u} we have $U_1=4$ and so by \eqref{Ugm} the number of primitive solutions of \eqref{secondinequality} with $0<\min(x,y) \leq Y_S$ is
$$\leq 8(12s+z(r)+\log h^{1/r}).$$
Taking into account solutions with $-Y_S\leq \min(x,y)<0$ and the solution $(1,0)$ we get
\be\label{sma}
P_{sma}(Y_S)\leq 1+16\big(12s+z(r)+\log h^{1/r}\big).
\ee}

\subsection{ Medium solutions}

For dealing with medium solutions, we invoke the following result which Mueller and Schmidt derives from \cite[Lemmas 15 and 16]{Mue-Sch} (see the comments in the beginning of \cite[Section 9]{Mue-Sch}).
\begin{lemma}
Let $(x,y)$ be a solution of \eqref{inequality} with 
$\min(|x|,|y|)\geq Y_S.$ Then  either there exists a root $\alpha$ of $F(x,1)$ satisfying
$$
\bigg|\alpha-\frac{x}{y}\bigg|<\frac{1}{H^{(1/s)-(1/r)}}\bigg(\frac{(rs)^{2s}(4e^3s)^rh}{|y|^r}\bigg)^{1/s}
$$
or there exists a root $\alpha^{-1}$ of $F(1,y)$  satisfying
$$
\bigg|\alpha^{-1}-\frac{y}{x}\bigg|<\frac{1}{H^{(1/s)-(1/r)}}\bigg(\frac{(rs)^{2s}(4e^3s)^rh}{|x|^r}\bigg)^{1/s}.
$$
\end{lemma}
Now we use \eqref{alphal1} to obtain
\begin{lemma}\label{17}
There is a set $S_1$ of roots of $F(x,1)$ and a set $S^*_1$ of roots of $F(1,y)$ both with cardinalities $\leq T(F)$ such that any solution of \eqref{inequality} with
$$\min(|x|,|y|)\geq Y_S$$
either has
\be\label{alpha}
\left|\alpha-\frac{x}{y}\right|\leq \frac{R(rs)^2}{H(F)^{\frac{1}{s}-\frac{1}{r}}}\bigg(\frac{(4e^3s)^rh}{|y|^r}\bigg)^{1/s}
\ee
for some $\alpha\in S_1$ or has
\be\label{alpha*}
\bigg|\alpha^*-\frac{x}{y}\bigg|\leq \frac{R(rs)^2}{H(F)^{\frac{1}{s}-\frac{1}{r}}}\bigg(\frac{(4e^3s)^rh}{|x|^r}\bigg)^{1/s}
\ee
for some $\alpha^*\in S^*_1.$
\end{lemma}
By the above lemma, it is enough to estimate the number of primitive solutions $(x,y)$ satisfying \eqref{alpha} for some $\alpha\in S_1$ or \eqref{alpha*} for some $\alpha*\in S^*_1$ with
\be\label{YSYL}
Y_S\leq y\leq Y_L\ {\textrm or}\ Y_S\leq x\leq Y_L,
\ee
respectively.
We consider the case when \eqref{alpha} is satisfied. The other case is similar. 
Our proof follows that of \cite[Lemma 5.1]{Sara-Div1}. 
Let $\alpha\in S_1.$ Let $V=\{(x_1,y_1),\ldots, (x_\nu,y_\nu)\} $ be the set of all solutions of 
\eqref{alpha} with $\gcd(x_i,y_i)=1$ and
$$Y_S\leq y_1\leq \cdots\leq y_\nu\leq Y_L.$$
Suppose that $\nu\geq 2.$ Then
\begin{align*}
\frac{1}{y_iy_{i+1}}\leq &\bigg|\frac{x_i}{y_i}-\frac{x_{i+1}}{y_{i+1}}\bigg|\leq 
\bigg|\alpha-\frac{x_i}{y_i}\bigg|+\bigg|\alpha-\frac{x_{i+1}}{y_{i+1}}\bigg|\\
\leq & \frac{K}{2y_i^{r/s}}+\frac{K}{2y_{i+1}^{r/s}}\leq \frac{K}{y_i^{r/s}}
\end{align*}
where
$$K=2R(rs)^2(4e^3s)^{r/s}h^{1/s}H(F)^{(1/r)-(1/s)}.$$
Thus we have 
$$y_{i+1}\geq K^{-1}y_i^{(r/s)-1}.$$
We now apply \cite[Lemma 2.1]{Sara-Div1} with
$$``T(x_i,y_i)=y_i, \beta=\frac{1}{K},\gamma=\frac{r-s}{s},A_1=Y_S,  B_1=Y_L$$
and
$$ \kappa=1\ {\rm or}\ 2\ {\rm according \ as}\ \beta\leq 1\
{\rm or}\ >1.''$$
Note that 
$$``\gamma"=\frac{r}{s}-1\geq 2$$
since $r\geq 3s.$ Then we have
\be\label{nu}
\nu\leq 1+\frac{1}{\log\big(\frac{r}{s}-1\big)}\frac{\log Y_L}{\log \big(\frac{Y_S}{K^{1/\kappa(\gamma-1)}}\big)}.
\ee
We shall now estimate the quantities involved in the above bound for $\nu.$
First
$$\frac {Y_S}{K^{\frac{1}{\kappa(\gamma-1)}}}\geq \frac{R^{\frac{s}{r-2s}}H(F)^{\frac{r-s}{2r(r-2s)}}}{(rs)^{\frac{2s}{r-2s}}}\geq 
R^{\frac{s}{2(r-2s)}}H(F)^{\frac{r-s}{2r(r-2s)}}.$$
We use the following estimates for $r\geq 3$ with the given choices of $(a,b).$  
$$\log R=z(r)+5;  \ \lambda\leq \frac{(1.64) \sqrt{r}}{(1-b)};\ r-\lambda\geq (.03)r.$$
Since $b\leq (.54),$ we get
$$\lambda\leq 3.566\sqrt{r}; \frac{1}{r-\lambda}\leq \frac{34}{r};  \frac{\lambda}{r-\lambda}\leq \frac{119}{\sqrt{r}}.$$
Using \eqref{Ma2}, $a\geq (.1)$ and $z(r)/r\leq 164.38$ (maximum occurring at $r=13$), 
we get
\begin{align*}
&\log Y_L \\
&\leq
\frac{1}{r-\lambda}\bigg((r+1)\log 2+\frac{r}{2}\log(r(r+1))+r\log H(F)+z(r)+3+\log h \bigg)+\\
&\frac{\lambda}{r-\lambda}\bigg(\log 4+\frac{1}{2a^2}\log(r+1)+\frac{1}{a^2}\log H(F)+\frac{r}{2a^2} \bigg)\\
&\leq  34\bigg(\big(1+\frac{1}{r}\big)\log 2+\frac{1}{2} \log \big(1+\frac{1}{r}\big)+\frac{5}{r}+\frac{z(r)}{r}+\log r+\log H(F)+\log h^{1/r} \bigg)+\\
&+119\bigg(\frac{\log 4}{\sqrt{r}}+\frac{50\log(r+1)}{\sqrt{r}}+\frac{100}{\sqrt{r}}\log H(F)+50\sqrt{r} \bigg)\\
&\leq 34\bigg(168+\log r+\log H(F)+\log h^{1/r}\bigg)+119\bigg(41+58\log H(F)+50\sqrt{r}\bigg)\\
& \leq 10591+34 \log r+5950 \sqrt{r}+6936 \log H(F)+34\log h^{1/r}.
\end{align*}
and
$$\log \bigg(\frac {Y_S}{K^{\frac{1}{\kappa(\gamma-1)}}}\bigg)\geq 
\frac{1}{r(r-2s)}\left(34rs\log^3r+\frac{r-s}{2}\log H(F)\right).$$
As $r\geq 3s,$ we get
\begin{align*}
\frac{\log Y_L}{\log \bigg(\frac {Y_S}{K^{\frac{1}{\kappa(\gamma-1)}}}\bigg)}
\leq & r^2\bigg(\frac{(10591+34 \log r+5950 \sqrt{r})}{34r\log^3r}+\frac{13872}{r-s}+\frac{68\log h^{1/r}}{(r-s)\log H(F)}\bigg)\\
\leq &r^2\bigg(7091+\frac{68\log h^{1/r}}{(r-s)\log H(F)}\bigg).
\end{align*}
So
$$\log \bigg(\log Y_L/\log \bigg(\frac {Y_S}{K^{\frac{1}{\kappa(\gamma-1)}}}\bigg)\bigg)\leq 9+2\log r+\log\bigg(1+\frac{\log h^{1/(r(r-s))}}{\log H(F)}\bigg).$$
Therefore by \eqref{nu}, we get
$$\nu\leq 1+\frac{9+2\log r+\log \bigg(1+\frac{\log h^{1/(r(r-s))}}{\log H(F)}\bigg)}{\log \gamma}.$$
Suppose $r\leq 3s^3.$ Then using $\gamma \geq 2,$ we find
$$\nu\leq 18+9\log^*s+(1.5) \log \bigg(1+\frac{\log h^{1/r(r-s)}}{\log H(F)}\bigg).$$
Since $|S_1|\leq 6s,$ we see that the number of primitive solutions satisfying \eqref{alpha} is 
$$\leq 162 s\log^* s+9s\log\bigg(1+\frac{\log h^{1/r(r-s)}}{\log H(F)}\bigg).$$
Suppose $r>3s^3.$ Then $r\geq 4$ and $s< (r/3)^{1/3}$ giving $\gamma=\frac{r}{s}-1\geq (.961)r^{2/3}.$ Hence $\log \gamma\geq (.62) \log r$ so that
$$\nu\leq 15+\frac{(1.62)}{\log r} \log\bigg(1+\frac{\log h^{1/r(r-s)}}{\log H(F)}\bigg).$$
This gives the number of primitive solutions satisfying \eqref{alpha} as 
$$\leq 90 s+8s\log\bigg(1+\frac{\log h^{1/r(r-s)}}{\log H(F)}\bigg)$$
since $r\geq 4.$
Similar estimates hold when \eqref{alpha*} is satisfied. Thus 
the number of primitive medium solutions satisfying \eqref{YSYL} is given by
$$
\begin{cases}
324 s\log^*s+18s\log\bigg(1+\frac{\log h^{1/r(r-s)}}{\log H(F)}\bigg)\ \textrm{if} \ r\leq 3s^3,\\
180s+16s\log\bigg(1+\frac{\log h^{1/r(r-s)}}{\log H(F)}\bigg)\ \textrm{if}\ r>3s^3.
\end{cases}
$$
Counting negative solutions as well, we see that
\be\label{medium}
P_{med}(Y_S,Y_L)\leq 648s(\log^*s)^\epsilon+36s\log \bigg(1+\frac{\log h^{1/r(r-s)}}{\log H(F)}\bigg).
\ee

\subsection{Final Assault for Theorem \ref{Main3}}
{\color{black}Let
$$G(m)=1+324s+648s (\log^*s)^\epsilon+16z(r)+16\log m^{1/r}+36s\log\bigg(1+\frac{\log m^{1/r}}{\log H(F)}\bigg).$$
Putting together the estimates from \eqref{lar}, \eqref{sma} and \eqref{medium}, we get
$$P'(F,h)\leq G(h).$$
Now Theorem \ref{Main3} follows from \eqref{Gm} and \eqref{GNFh}.
\qed}
\section{Correction in the paper of Bengoechea \cite{Ben} }
In this final section we point out an error in the following result of Bengoechea \cite[Theorems 1.1, 1.2, Corollary 1.3]{Ben}.
\vskip 2mm
\noindent
{\bf Theorem. B} (Bengoechea). {\it Suppose $F$ has at most $s+1$ non-zero terms.
\vskip 2mm
\noindent
(i) Suppose 
\be\label{8}
|D(F)|>(r(r-1))^{8r(r-1)}.
\ee
Then
$$N(F,h)\ll sh^{2/r}.$$
\vskip 2mm
\noindent
(ii) Suppose $r\geq 3s.$ Then
$$N(F,h)\ll \bigg((c(s)(1+\log h^{1/r})+\log^3r\bigg)h^{2/r}|D(F)|^{-\frac{1}{r(r-1)}}$$
where
\be\label{cs}
c(s)=\begin{cases}
s\ \textrm{ if } \ r\geq s^4,\\
s\log s\ \textrm{ if } \ 9s^2\leq r<s^4,\\
s\log s\bigg(1+\frac{s}{\log H(F)}\bigg)\ \textrm{ if }\ r<9s^2.
\end{cases}
\ee
\vskip 2mm
\noindent
(iii) Suppose $r\geq 3s$ and 
$$
|D(F)|^{\frac{1}{2.5(r-1)}}\geq h,
$$ 
then 
$$N(F,h)\ll c(s)+\log^3r.$$
}
Needless to say, her work depends on the methods in \cite{Sch} and \cite{Mue-Sch}. She combines the two hypotheses of forms having $s+1$ non-zero coefficients and forms having minimal Mahler measure judiciously to get new bounds especially for {\it small} 
solutions. 
\vskip 2mm
\noindent
We consider the estimates on pages 1233 and 1234 in  \cite{Ben}. 
In the notation of \cite{Ben}, we have $n=r,h=m$ and $\tilde P=P'.$
In \cite[equation (43)]{Ben}, the author chooses the smallest prime $p$ that satisfies
$$p\geq e^{400} h^{2/r}|D(F)|^{-1/r(r-1)}.$$
Then by Bertrand's postulate, she derives in \cite[equation (44)]{Ben} that
\be\label{400}
p<2e^{400}h^{2/r}|D(F)|^{-1/r(r-1)}.
\ee
In \cite[p. 1234]{Ben}, she estimates $ P'(F_{A_j},h)$ as
$$P'(F_{A_j},h)\ll r+s\log\log h^{1/r}$$
and so deduces in \cite[equation (46)]{Ben} that
\be\label{p1}
P'(F,h)\ll (p+1)(r+s\log\log h^{1/r})\ll \left(r+s\log\log h^{1/r}\right)h^{2/r}|D(F)|^{-1/r(r-1)}
\ee
as seen in \eqref{pinequality}. After this she considers two cases according as $|D(F)|^{1/r(r-1)}>\log\log h^{1/r}$ or
$|D(F)|^{1/r(r-1)}\leq \log\log h^{1/r}$ to conclude that $P'(F,h)\ll sh^{2/r}.$ The latter case depends on a theorem of Thue and \eqref{8}. The estimate for $N(F,h)$ now follows by a standard retracing procedure. See \eqref{NPP} also.

Note that \eqref{400} is true only when 
$$e^{400}h^{2/r}|D(F)|^{-1/r(r-1)}> 1.$$ 
In the complementary case one needs to take $p=2$ and then \eqref{400} (and hence \cite[equation (44)]{Ben} )is not valid. For instance, consider $(i)$ of Theorem B and put
$$D_0:=e^{400}h^{2/r}/(r(r-1))^8.$$
Then by \eqref{8}, 
$$e^{400}h^{2/r}|D(F)|^{-1/r(r-1)}\leq D_0\leq 1$$ 
can very well occur for large $r$ and small $h.$ 
In order to  avoid this discrepancy, we may replace \eqref{400} with
$$p\leq\max\left(2,2e^{400}h^{2/r}|D(F)|^{-1/r(r-1)}\right)<2\left(1+e^{400}h^{2/r}|D(F)|^{-1/r(r-1)}\right).$$
This would then affect \eqref{p1} (and hence \cite[equation (46)]{Ben}) as
$$ P' (F,h)\ll \left(1+\frac{h^{2/r}}{|D(F)|^{1/r(r-1)}}\right)\left(r+s\log\log h^{1/r}\right).$$

Similar discrepancy is found in p. 1239 (see equation (60) and the last line) which results in a correction on p.1240, l. 11 of \cite{Ben} 
as
$$\ll \left(1+\frac{h^{2/r}}{|D(F)|^{1/r(r-1)}}\right) \left(c(s)(1+\log h^{1/r})+\log^3r\right).$$ 
This affects the results $(ii)$ and $(iii)$ of Theorem B.
Essentially the error occurs in the application of Bertrand's result and it completely changes the quality of the result! The above discussion leads to the following correction of the result of Bengoechea.
\vskip 2mm
\noindent
{\bf Theorem} (Correct version). {\it Suppose $F$ has at most $s+1$ non-zero terms.
\vskip 2mm
\noindent
(i) Suppose 
$$|D(F)|>(r(r-1))^{8r(r-1)}.$$
Then
$$N(F,h)\ll\left( r+s\log\log h^{1/r}\right) h^{2/r}.$$
\vskip 2mm
\noindent
(ii) Suppose $r\geq 3s.$ Then
$$N(F,h)\ll  \left(1+\frac{h^{2/r}}{|D(F)|^{1/r(r-1)}}\right) \left(c(s)(1+\log h^{1/r})+\log^3r\right)h^{2/r}$$
where $c(s)$ is given by \eqref{cs}.
\vskip 2mm
\noindent
(iii) Suppose $r\geq 3s$ and 
$$|D(F)|^{\frac{1}{2.5(r-1)}}\geq h,$$  
then 
$$N(F,h)\ll \left(c(s)(1+\log h^{1/r})+\log^3r\right)h^{2/r}.$$
}

In fact, the assertion in $(iii)$ above is true if $|D(F)|^{\frac{1}{2(r-1)}}\geq h,$ as can be seen from $(ii).$  

\section*{Acknowledgment}
Divyum acknowledges the support of the DST--SERB SRG Grant (SRG/2021/000773) and the OPERA award of BITS Pilani.


\end{document}